\documentclass[a4paper,12pt,reqno]{amsart}
\usepackage{amssymb}
\usepackage[dvips]{graphicx}
\usepackage{hyperref}
\usepackage{inputenc}

\input epsf.tex
\newdimen\xsize
\newdimen\oldbaselineskip
\newdimen\oldlineskiplimit
\def\restorelineskip{\baselineskip=\oldbaselineskip%
\lineskiplimit=\oldlineskiplimit}
\def\putm[#1][#2]#3{% #1 X-axis  #2 Y-axis(downwords)  #3 $"mathobject"$
\hbox{\vbox to 0pt{\parindent=0pt%
\vskip#2\xsize\hbox to0pt{\hskip#1\xsize $#3$\hss}\vss}}}%
\long\def\Line#1{\hbox to \hsize{#1}}
\def\putt[#1][#2]#3{% #1 X-axis  #2 Y-axis(downwords)  #3 "text"
\vbox to 0pt{\noindent\hskip#1\xsize\lower#2\xsize%
\vtop{\restorelineskip#3}\vss}}

\makeatletter
\def\xbig[#1]#2{{\hbox{$\m@th\left#2\vbox to#1\xsize{}%
\right.\n@space$}}}
\makeatother
\def\xlar[#1]#2{%
\smash{\mathop{ \hbox to #1\xsize{\leftarrowfill}}\limits^{#2}}}
\def\xrar[#1]#2{%
\smash{\mathop{ \hbox to #1\xsize{\rightarrowfill}}\limits^{#2}}}
\def\xline[#1]{\hbox to #1\xsize{\leaders\hrule\hfill}}

\DeclareFontFamily{U}{rsf}{\skewchar\font'177}%
\DeclareFontShape{U}{rsf}{m}{n}{<-6>rsfs5<6-8>rsfs7<8->rsfs10}{}%
\DeclareFontShape{U}{rsf}{b}{n}{<-6>rsfs5<6-8>rsfs7<8->rsfs10}{}%
\DeclareMathAlphabet\RSFS{U}{rsf}{m}{n}
\SetMathAlphabet\RSFS{bold}{U}{rsf}{b}{n}
\def\sf#1{{\mathsf{#1}}}%\def\sf{\fam\sffam}

\def\slsf{\slshape \sffamily }

\def\msmall#1{\mathchoice{\hbox{\small$\displaystyle {#1}$}}{#1}{#1}{#1}}

\def\qq{{\mathbb Q}}

\def\s{{\mathbf S}}

\def\b{{\mathbf B}}
\def\cc{{\mathbb C}}

\def\rr{{\mathbb R}}

\def\nn{{\mathbb N}}

\def\zz{{\mathbb Z}}

\def\a{\sf{a}}
\def\b{\sf{b}}
\def\adyn{\sf{1}}
\def\arg{{\sf{Arg}}}

\def\sym{\sf{sym}}

\def\C{\sf{C}}

\def\cos{\sf{cos}\,}
\def\const{\sf{const}}

\def\dist{\sf{dist}\,}
\def\e{\sf{e}}

\def\el{\sf{l}}

\def\s{\sf{s}}
\def\sym{\sf{Sym}}

\def\im{\sf{Im}\,}

\def\re{\sf{Re}\,}
\def\k{\sf{K}}

\def\lim{\mathop{\sf{lim}}}
\def\limsup{\mathop{\sf{lim\,sup}}}
\def\Ln{\sf{Ln}}
\def\log{\sf{log}\,}

\def\n{{\mathrm{n}}}

\def\max{\sf{max}}
\def\min{\sf{min}}

\def\p{\sf{p}}

\def\res{\sf{Res}}
\def\r{{\mathrm{r}}}
\def\x{{\mathrm{x}}}
\def\\xi{\mathrm{\xi}}
\def\ro{{\mathrm{\rho}}}

\def\sin{\sf{sin}\,}

\def\sup{\sf{sup}\,}

\def\w{{\mathrm{w}}}

\def\z{\sf{z}}

\def\eps{\varepsilon}

\def\<{\langle}\let\la=\<
\def\>{\rangle}\let\ra=\>
  
\def\comp{\Subset}
\def\d{\partial}

\def\ddef{\mathrel{{=}\raise0.3pt\hbox{:}}}
\def\deff{\mathrel{\raise0.3pt\hbox{\rm:}{=}}}

\def\fraction#1/#2{\mathchoice{{\msmall{ #1\over#2}}}%
{{ #1\over #2 }}{{#1/#2}}{{#1/#2}}}
\def\norm#1{\left\Vert{#1}\right\Vert}

\def\le{\leqslant}

\def\emptyset{\varnothing}

\def\longpoints{\leaders\hbox to 0.5em{\hss.\hss}\hfill \hskip0pt}
\def\stateskip{\smallskip}
\def\state#1. {\stateskip\noindent{\bf#1. }} %\medskip
\def\statep#1. {\stateskip\noindent{\bf#1 }} %\medskip
\def\proof{\state Proof. \2}

\def\Chi{\raise 2pt\hbox{$\chi$}}

\def\ie{\hskip1pt plus1pt{\sl i.e.\/,\ \hskip1pt plus1pt}}

\def\sli{{\sl i)} } 
\def\slii{{\sl i$\!$i)} } 
\def\sliii{{\sl i$\!$i$\!$i)} }

\def\Chi{\raise 2pt\hbox{$\chi$}}

\let\phI=\phi\let\phi=\varphi\let\varphi=\phI
\let\cal=\mathcal

\def\calc{{\cal C}}
\def\cald{{\cal D}}

\def\call{{\cal L}}
\def\calm{{\cal M}}

\def\calp{{\cal P}}

\def\calr{{\cal R}}

\def\calz{{\cal Z}}

\def\eps{\varepsilon}

\def\comp{\Subset}
\def\d{\partial}

\def\1{{1\mkern-5mu{\rom l}}}

\def\ge{\geqslant}

\def\fraction#1/#2{\mathchoice{{\msmall{ #1\over#2}}}%
{{ #1\over #2 }}{{#1/#2}}{{#1/#2}}}

\def\le{\leqslant}

\def\emptyset{\varnothing}
\newcommand{\2}{\thinspace}

\def\qed{\ \ \hfill\hbox to .1pt{}\hfill\hbox to .1pt{}\hfill $\square$\par}

\def\comment#1\endcomment{}

 \belowdisplayskip=\abovedisplayskip

\def\lineeqqno(#1){\hfill\llap{\vbox to 10pt%
{\vss\begin{align} \eqqno(#1)\end{align}\vss}}\vskip1pt}

\advance\headheight 1.2pt

\def\ShowwLLabel#1{}

\def\thechpt{\Roman{chpt}}
\def\newchapt[#1]#2{%
\refstepcounter{chpt}\setcounter{subsection}{0}%
\setcounter{thm}{0}\setcounter{defi}{0}%
\setcounter{rema}{0}\setcounter{exrc}{0}%
\renewcommand{\thesubsection}{\thechpt.\arabic{subsection}}%
\newpage
\section*{\begin{center}\huge \bf Chapter \thechpt\\
#2 \end{center}}\label{#1}%
\ \smallskip%
\addcontentsline{toc}{chpt}{Chapter \thechpt. #2}%
\markboth{Chapter \thepart}{#2}%
}
\def\newsect[#1]#2{\refstepcounter{section}\setcounter{equation}{0}%
\renewcommand{\thesubsection}{\arabic{section}.\arabic{subsection}}%
\section*{\arabic{section}.
#2}\vspace{-20pt}\label{#1}\vspace{20pt}%
\markboth{Section \arabic{section}}{#2}}

\def\newlect[#1]#2{\refstepcounter{section}%
\renewcommand{\thesubsection}{\arabic{section}.\arabic{subsection}}%
\section*{Lecture \arabic{section}\\
#2}\label{#1}%
\markboth{Lecture \arabic{section}}{#2}}

\def\newprg[#1]#2{\refstepcounter{subsection}%
\subsection*{{\thesubsection.\ #2}} \label{#1}%
}

\def\newprgnew[#1]#2{\refstepcounter{subsection}\setcounter{equation}{0}%
\subsection*{{\mdseries\slshape \sffamily \thesubsection.\ #2}} \label{#1}%
\addtocontents{toc}{{\bf\thesubsection.} #2,~p.\pageref{#1}.\ \ }%
}

\def\newappx[#1]#2{\refstepcounter{appx}\setcounter{section}{0}%
\renewcommand{\thesubsection}{A\arabic{appx}.\arabic{subsection}}%
\section*{Appendix \arabic{appx} #2}%\Roman{appx}??
\label{#1}%
\markboth{Appendix A\arabic{appx}}{#2}
}

\newtheorem{thm}{Theorem}[section]
   \def\newthm#1{\begin{thm}\label{#1}}

\newtheorem{nnthm}{Theorem}
   \def\newthm#1{\begin{nnthm}\label{#1}}

\newtheorem{lem}{Lemma}[section]
   \def\newlemma#1{\begin{lem} \label{#1}}

\newtheorem{prop}{Proposition}[section]
   \def\newprop#1{\begin{prop}\label{#1}}

\newtheorem{nnprop}{Proposition}
   \def\newprop#1{\begin{nnprop}\label{#1}}

\newtheorem{corol}{Corollary}[section]
   \def\newcorol#1{\begin{corol} \label{#1}}

\newtheorem{nncorol}{Corollary}
   \def\newcorol#1{\begin{nncorol} \label{#1}}

\newtheorem{defi}{Definition}[section]
   \def\newdefi#1{\begin{defi} \label{#1}\rm }

\newtheorem{nndefi}{Definition}
   \def\newdefi#1{\begin{defi} \label{#1}\rm }

\newtheorem{exmp}{Example}[section]
   \def\newexmp#1{\begin{exmp} \label{#1}\rm }

\newtheorem{nnexmp}{Example}
   \def\newexmp#1{\begin{nnexmp} \label{#1}\rm }

\newtheorem{exrc}{Exercise}
   \def\newexrc#1{\begin{exrc} \label{#1}\rm }

\newtheorem{rema}{Remark}[section]
   \def\newrema#1{\begin{rema} \label{#1}\rm }

\newtheorem{nnrema}{Remark}
   \def\newthm#1{\begin{nnrema}\label{#1}}

\def\eqqno(#1){\label{(#1)}}
\def\eqqref(#1){(\ref{(#1)})}

\DeclareMathOperator*{\Sup}{sup}
\makeindex

\title{Riemann surface of the Riemann zeta function}
\author{S. Ivashkovich}
\date{\today}

\address{
Universit\'e de Lille-1, UFR de Math\'ematiques, 59655 Villeneuve
d'Ascq, France} \email{ivachkov@math.univ-lille1.fr}

\begin{document}
\begin{abstract}
In this paper we treat the classical Riemann zeta function as a function of three variables: one is 
the usual complex $\adyn$-dimensional, customly denoted as $s$, another two are complex infinite dimensional,
we denote it as $\b = \{b_n\}_{n=1}^{\infty}$ and $\z =\{z_n\}_{n=1}^{\infty}$. When $\b = \{1\}_{n=1}^{\infty}$
and $\z = \{\frac{1}{n}\}_{n=1}^{\infty}$ one gets the usual Riemann zeta function. Our goal in this paper  
is to study the meromorphic continuation of $\zeta (\b , \z ,s)$ as a function of the triple $(\a , \z , s)$.
\end{abstract}

\maketitle

\setcounter{tocdepth}{1}
\tableofcontents

%\addtocontents{toc}{

%\bigskip}
%\addtocontents{toc}{ \leftskip=0mm \rightskip=0mm}

\newsect[INT]{Introduction.}

\newprg[INT.e-r]{Euler-Riemann zeta function}
For a real $\sigma$ let $H_{\sigma} \deff \{s \in\cc:\re s > \sigma\}$ be the halph-plane. For a complex variable 
$s\in H_1$ the classical Riemann or, better Euler-Riemann zeta function is defined as 
\begin{equation}
\eqqno(zeta-f1)
\zeta (s) = \sum_{n=1}^{\infty}\frac{1}{n^s} = \sum_{n=1}^{\infty}\left(\frac{1}{n}\right)^s = 
\zeta \left(\frac{\adyn}{\n}, s\right).
\end{equation}
In the right hand side of \eqqref(zeta-f1) we interpret $\frac{\adyn}{\n}\deff \{\frac{1}{n}\}_{n=1}^{\infty}$ 
as the distinguished element of the following complex {\slsf sequence space } 
\begin{equation}
\eqqno(1/n-seq)
\el_{\frac{\adyn}{\n}} \deff \{\z = \{z_n\}_{n=1}^{\infty}: \norm{\z}_{\frac{\adyn}{\n}} \deff
\sup_{\substack{n\in \nn}} \left(n|z_n|\right)<\infty\}.
\end{equation}
We shall also use the notation $\n^{-1}$ for $\frac{\adyn}{\n}$ as well as $\el_{\n^{-1}}$ for 
$\el_{\frac{\adyn}{\n}}$. More generally:

\begin{nndefi}
\label{el-r1}
Let $\r = \{r_k\}_{k\in \nn}$ be a sequence of positive reals. The following space of complex sequences
\begin{equation}
\eqqno(el-r2)
\el_{\r} \deff \left\{\z=\{z_n\}_{n=1}^{\infty}: \norm{\z}_{\r}\deff \sup_{\substack{n\in \nn}}
\left|\frac{z_n}{r_n}\right| <\infty \right\}
\end{equation}
we shall call the $\r$-sequence space or, simply an $\r $-space.
\end{nndefi}

\noindent $\el_{\r}$ is a Banach  space with respect to the norm $\norm{\cdot}_{\r}$. In fact it is isometrically 
isomorphic to $\el^{\infty}$. Actually $\el^{\infty}$ appears to be $\el_{\adyn}$ according to this definition,
here $\adyn\deff \{1\}_{n=1}^{\infty}$.  By $B^{\infty}_{\r}(\z^0,\beta)$ denote the open ball of radius $\beta >0$ in 
the Banach space $\el_{\r}$ centered at $\z^0$. Notice that for $\z =\{z_n\}\in B^{\infty}_{\n^{-1}}(\n^{-1}, \adyn)$ we have 
that $\re z_n>0$ and therefore we can fix that $\arg z_n\in ]-\pi /2,\pi /2[$.

\smallskip For given $\b=\{b_n\}_{n=1}^{\infty}\in \el^{\infty}$, $\z=\{z_n\}_{n=1}^{\infty}\in B^{\infty}_{\n^{-1}}(\n^{-1}, \adyn)$ 
and $s\in H_1$ set (with arguments of $z_n$ taken as above):
\begin{equation}
\eqqno(zeta-f2)
\zeta^{ER}(\b, \z, s) = \sum_{n=1}^{\infty}b_nz_n^s
\end{equation}
and call it still the Euler-Riemann zeta function. It is not difficult to observe, see Lemma \ref{zeta-def},
that $\zeta^{ER}$ is holomorphic as the function of the triple $(\b,\z,s)$ on $\el^{\infty}\times B^{\infty}_{\n^{-1}}(\n^{-1}, \adyn)
\times H_1$. Variable $\z$ we shall call a {\slsf vectorial variable} of $\zeta^{ER}$ and $\b$ a (complex) {\slsf parameter}. 
$\zeta^{ER}(\b,\z,s)$ will be denoted simply as $\zeta (\b,\z,s)$ whenever this will not lead to a confusion. Notice that 
$\zeta (\b,\z,s)$ is certainly {\slsf not} well defined on the coordinate cross $\C^{\infty}_{\n^{-1}} = 
\{\z = \{z_n\}_{n=1}^{\infty}$ $\in \el_{\n^{-1}}: z_n=0 \text{ for some } n\}$ at least for a general 
$\b\in \el^{\infty}$. The closure of $\C^{\infty}_{\n^{-1}}$ in $\el_{\n^{-1}}$ is 
\[
\bar\C^{\infty}_{\n^{-1}} = \left\{\z\in \el_{\n^{-1}}:\inf_{\substack{n\in \nn}}
n|z_n|=0\right\}, 
\]
see Lemma \ref{c-cross}. We prove in Lemma \ref{l2-cont} that $\zeta$ can be 
continued to a (multivalued) analytic function on $\el^{\infty}\times \left(\el_{\n^{-1}}\setminus \bar\C_{\n^{-1}}\right)
\times H_1$, again  as a function of the triple $(\b , \z , s)$. 

\begin{nnrema} \rm
\label{b-arg0}
One should be a bit carefull when speaking about multivalued analyticity and universal covers in infinite dimensions.
In our case the space $\el_{\n^{-1}}\setminus \bar C^{\infty}_{\n^{-1}}$ admits the universal 
cover since it is semi-locally simply connected, see Corollary 14 in \cite{Sp}. 
\end{nnrema}

Denote by  $\calm$ the maximal set over 
$\el^{\infty}\times \left(\el_{\n^{-1}}\setminus \bar\C^{\infty}_{\n^{-1}}\right)$  which consists of 
such $(\b,\z)$ that $\zeta (\b,\z, \cdot)$ can be analytically continued to a meromorphic function on 
$\cc$. By the classical theorem of Riemann $(\adyn , \n^{-1})\in \calm$, here $\adyn$ 
and $\n^{-1}$ are sequences defined as above. Furthermore denote by $\calz$ the maximal set over 
$\el^{\infty}\times \left(\el_{\n^{-1}}\setminus \bar\C^{\infty}_{\n^{-1}}\right)\times \cc$ to which 
$\zeta (\b,\z,s)$ can be analytically continued as a meromorphic function, \ie $\calz$ is the ``Riemann surface'' 
of $\zeta$. Notice that $\calz \supset \calm\times \cc$.
In this paper we are interested in the following:

\smallskip\noindent
{\slsf Question: \it What can be said about the structure of $\calm$ and $\calz$?}

\smallskip  Our first result states that $\calm$ is ``quite thick'' and, moreover, possesses a certain
Banach analytic structure. Denote by $\el_{\e^{-\n}}$ 
the sequence space defined by the sequence $\e^{-\n} \deff \{e^{-n}\}_{n\in \nn}$, \ie 
\begin{equation}
\eqqno(e-n-seq)
\el_{\e^{-\n}} \deff \{ \{z_n\}_{n=1}^{\infty}: \norm{\z}_{\e^{-\n}} \deff 
\Sup_{\substack{n\in \nn}} \left(e^n|z_n|\right)<\infty\}.
\end{equation}
Notice that $\el_{\e^{-\n}}$ is obviously densily imbedded into  $\el_{\n^{-1}}$.
We prove the following:

%\newpage

\begin{nnthm}
\label{m-thick}
For every $(\b^0, \z^0)\in \calm$ the following holds: 

\sli $(\b^0 ,\z^0 + B_{\e^{-\n}}^{\infty}(0,\norm{\z^0}_{\n^{-1}})) \subset \calm$, 
or in other words 
$\{(\b^0 ,\z^0 +  B_{\e^{-\n}}^{\infty}(0,\norm{\z^0}_{\n^{-1}}))\} \times \cc \subset \calz$;

\slii the restriction $\zeta (\b^0 , \cdot , \cdot)$ of $\zeta$ to $\{\b^0\}\times \{\z^0+ 
B_{\e^{-\n}}^{\infty}(0,\norm{\z^0}_{\n^{-1}})\}\times \cc$ can  be analytically

\quad   continued as a multivalued meromorphic function to 
$\{\b^0\}\times\left(\{\z^0 + \el_{\e^{-\n}}\}\setminus \bar\C^{\infty}_{\n^{-1}}\right)\times \cc$.
\end{nnthm}

\begin{nnrema} \rm

\smallskip\noindent{\bf 1.}  First of all let's make clear that for $\z^0\in \el_{\n^{-1}}$ one has
\[
B^{\infty}_{\e^{-\n}}(0,\norm{\z^0}_{\n^{-1}}) = \{\z\in \el_{\e^{-\n}}: \norm{\z}_{\e^{-\n}} < 
\norm{\z^0}_{\n^{-1}}\}.
\]

\smallskip\noindent{\bf 2.} When writing $(\b^0 , \z^0)\in \calm$ we mean, in particular, that $(\b^0 , \z^0)
\in \el^{\infty}\times (\el_{\n^{-1}}\setminus \bar\C_{\n^{-1}})$. Therefore due to Lemma \ref{two-crosses} 
we have that $\{\z^0 + \el_{\e^{-\n}}\}\setminus \bar\C^{\infty}_{\n^{-1}} = \{\z^0 + \el_{\e^{-\n}}\}\setminus 
\C^{\infty}_{\n^{-1}}$ and, moreover, in the same lemma we prove that $\{\z^0 +\el_{\e^{-\n}}\}\cap 
\C^{\infty}_{\n^{-1}}$ is a countable union of hyperplanes such that every bounded subset of the affine 
plane $\{\z^0 +\el_{\e^{-\n}}\}$ meets only finitely many of them. Therefore an analytic continuation within 
its complement makes well sense. 

\smallskip\noindent{\bf 3.} Let $L_n \deff \{\z=\{z_m\}_{m=1}^{\infty}\in \el_{e^{-\n}}: z_n + z_n^0=0\}$ 
be one of  these hyperplanes. Then it is easy to see that the monodromy of $\zeta (\b^0,\cdot , s)$ around 
this hyperplane along a loop $\gamma_n$ in $z_n$-plane starting at point $z^0_n+z_n\not= 0$ and going in the 
positive direction is $b^0_n(z^0_n+z_n)^s\left(e^{2\pi i s} - 1\right)$. Since $\{\z^0+\el_{e^{-\n}}\}\cap 
\C^{\infty}_{\n^{-1}} = \bigcup_nL_n$ the fundamental group of $\{\z^0+\el_{e^{-\n}}\}\setminus \C^{\infty}_{\n^{-1}}$
is generated by such $\gamma_n$-s and is abelian. This shows that $\zeta (\b^0,\cdot , \cdot)$ lifts to the 
universal cover of $\{\b^0\}\times \left(\{\z^0+\el_{e^{-\n}}\}\setminus \C^{\infty}_{\n^{-1}}\right)\times \cc$
and separates points in the fiber (for a non-integer $s$). 

\smallskip\noindent{\bf 4.} Since $\{\z^0 + \el_{\e^{-\n}}\}\times \cc$ is a dense affine subspace 
of $\el_{\n^{-1}}\times \cc$ this statement means that the Riemann surface of the function $\zeta$ contains a huge 
quantity of domains spread over parallel dense affine subspaces of the total space $\el^{\infty}\times
\el_{\n^{-1}}\times \cc$. Indeed, apart from $(\b^0,\z^0)= (\adyn,\n^{-1})$ we know quite
a number of pairs $(\b^0 , \z^0)\in \calm $. If, for example, $\alpha =\sqrt{D}$ is a quadratic irrationality  
then $(\{n\alpha\}_{n=1}^{\infty}, \n^{-1})\in \calm$ due to results of Hecke and Hardy-Littlewood,
see \cite{He, HL}. More examples can be found in \cite{Pa}, \cite{KL} and many other places.

\smallskip\noindent{\bf 5.} Notice as well that the result of this theorem can be reformulated in terms
of the general Dirichlet series, see subsection \ref{VV.dir}.

\smallskip\noindent{\bf 6.} It should be noticed that the analytic continuation in this theorem preserves
the poles of $\zeta (\b^0,\z^0,s)$.
\end{nnrema}

\newprg[INT.ll]{Zeta function in the form of Lerch-Lipschitz}
In order to produce new pairs $(\b,\z)\in \calm $  let us make the 
following substitution to $\zeta$:
\begin{equation}
\eqqno(chg-var)
\begin{cases}
b_n = e^{2\pi i na_n}\cr
z_n = \frac{1}{n+\xi_n}.
\end{cases}
\end{equation}
We obtain the function $\zeta$ in the form of Lerch-Lipschitz
\begin{equation}
\eqqno(ll-f0)
\zeta^{LL}(\a, \z ,s) = \sum_{n=1}^{\infty}\frac{e^{2\pi in a_n}}{(n + z_n)^s}.  
\end{equation}
Variable $\xi$ we immeediately renamed as $\z$, but parameter $\a$ will distinguish the Lerch-Lipschits form
of the function zeta from that of Euler-Riemann and $\zeta^{LL}$ will be renamed back to $\zeta$. This function 
is certainly well defined and holomorphic for $\a = \{a_n\}\in \el^{\infty}_{\im^+}\deff \{\a\in \el^{\infty}:
\im a_n > 0\}$, $\z = \{z_n\}\in \el^{\infty}_{\re^+} \deff \{\z\in \el^{\infty}: \re z_n > 0\}$ and $s\in H_1$,
see Corollary \ref{l2-cont1}. Consider the following ``$k$-dimensional cut`` $\zeta_k$ of function $\zeta $:
%, first for all real values of $(\a, \z)$
\begin{equation}
\eqqno(ll-f)
\zeta_k(\a, \z, s) = \sum_{n=1}^{\infty}\frac{e^{2\pi in a_{[n-1]_k+1}}}{(n + z_{[n-1]_k+1})^s}
= \frac{e^{2\pi ia_1}}{(1+z_1)^s} + ... + \frac{e^{2\pi ika_k}}{(k+z_k)^s} + \frac{e^{2\pi i(k+1)a_1}}{(k+1+z_1)^s}
+ ...
\end{equation}
Here $\a=(a_1,...,a_k)\in \cc^k_{\im^+}$, $\z=(z_1,...,z_k)\in \cc^k_{\re^+}$ and $[m]_k$ is the remainder of 
the division of $m$ by $k$. I.e., $\zeta_k$ is, in fact, $\zeta$ itself but restricted to the finite dimensional 
subspace of {\slsf ''$k$-periodic sequences``} $\left\{\left(a_{[n-1]_k+1}, z_{[n-1]_k+1}\right)
\right\}_{n=1}^{\infty}$. Our second result is the following:

\begin{nnthm}
\label{ll-thm}
Function $\zeta_k$ admits an analytic continuation to a single-valued meromorphic function $\tilde\zeta_k$ on 
\begin{equation}
\tilde\calz_k \deff \widetilde{\cc}^{k,ab}_{\zz [1/k]}\times \widetilde{\cc}^{k,ab}_{\zz^-}\times \cc .
\end{equation}
%In addition $\tilde \zeta_k$ separates points of $\tilde\calz_k$ over $\calz_k\deff  \cc^k_{\zz [1/k]}\times 
%\cc^k_{\zz^-}\times \cc$.
\end{nnthm}

\noindent Here we denote by $\zz [1/k]$ the subring of $\qq$ which consists
from rationals of the form $\{l/k:l\in \zz\}$ for a fixed $k$. $\widetilde{\cc}^{k,ab}_{\zz[1/k]}$ stands for
the abelian cover of $\cc^k_{\zz[1/k]} \deff \left(\cc\setminus \zz \left[1/k\right]\right)^k$. 
$\widetilde{\cc}^{k,ab}_{\zz^-}$ is the abelian cover of $\cc^k_{\zz^-}\deff (\cc\setminus \zz^-_1)\times ... 
\times (\cc \setminus \zz^-_k)$, where $\zz^- = \zz_{<0}$ is the set of negative integers and $\zz^-_j\deff 
\{z\in \zz^-: [z]_k=-j\}$.
%If $\pi:\tilde\calz\to \calz$ is a regular cover we say that function $\chi$ separates points of $\tilde\calz$zz^-
%over $\calz$ if for every $z\in \calz$ the germs $\chi_{\tilde z}$ are pairwise dinstinct for all preimages 
%$\tilde z$ of $z$ under $\pi$. 
The statement of the theorem means that the monodromy of $\zeta_k$ along any loop generating the commutator 
of $\pi_1(\cc^k_{\zz [1/k]}\times \cc^k_{\zz^-})$ is trivial and that $\tilde\calz_k$ is the genuine Riemann surface 
of $\zeta_k$. Moreover, this theorem shows that for any $k$ the set $\calm$ is infinitely sheeted over the space 
of $k$-periodic sequences. This provides us a huge quantity of new pais $(\a, \z)\in \calm$.  
\begin{nnrema} \rm
Periods of $a$-s and $z$-s in Theorem \ref{ll-thm} can be taken different, say $k_1$ and $k_2$. Then setting 
$k=lcm(k_1,k_2)$ we can apply our theorem. 
\end{nnrema}

\newprg[INT.proof]{The structure of the paper and proofs}
We organize the material of this paper as follows. First, in section \ref{VV} we recall and prove few topological 
properties of sequence spaces in order to clarify the statement of Theorem \ref{m-thick}. After that we 
prove this theorem. The needed material on complex and pluricomplex analysis on sequence spaces is
moved to the Appendix. In particular, the following infinite dimensional version of  Hartogs-Siciak theorem  will 
be proved there (\ie in Appendix). Given complex sequence spaces $\el_{\r_1}$ and $\el_{\r_2}$.

\begin{nnthm}
\label{hart-sic}
Let a holomorphic function $f$ on $B^{\infty}_{\r_1}(0,1)\times B^{\infty}_{\r_2}(0,1)$ be such that for every
$\z$ in some plurithick subset $E$ of $B^{\infty}_{\r_1}(0,1)$ the restriction $f(\z,\cdot)$ can be  continued to a 
holomorphic function on $\el_{\r_2}$. Then $f$ can be continued to a holomorphic function
on $B^{\infty}_{\r_1}(0,1)\times \el_{\r_2}$. 
\end{nnthm}

\smallskip Second, we prove Theorem \ref{ll-thm} in section \ref{LL}, it doesn't requires any infinite dimensional 
analysis and is accessible right away, \ie withought section \ref{VV}. Our proof follows the main lines of the original 
approach of Riemann \cite{Rm} and then of Lerch \cite{Le} and, finally of Lagarias and Li in \cite{LL}, where the case 
$k=1$ of Theorem \ref{ll-thm} was proved. The particularity of our approach is again the use of (finite dimensional)
Hartogs-Siciak theorem which makes the exposition considerably shorter and accessible when $k>1$.

\newsect[VV]{Vectorial variable in the Riemann zeta function}

We start this section recalling the notion of a sequence space and pointing out few facts about their
vectorial topology. They are absolutely needed for the very understanding of the statement of Theorem 
\ref{m-thick} from the Introduction as well as for its proof. After that we turn to the Riemann zeta 
function and to the proof of the mentioned theorem. Needed facts from complex and pluricomplex analysis 
on sequence spaces are placed to the Appendix.

\newprg[VV.sec]{Vectorial topology on sequence spaces}

Throughout this paper we consider various complex vector spaces of infinite sequences of complex numbers 
$\z=\{ z_k\}$, they are customly called the {\slsf sequence spaces}. If, for some reasons, we shall need 
to consider real sequences in the complex sequence space $\el_{\r}$ then the corresponding real vector space
will be denoted as $\el_{\r ,\rr}$. The standard basis of {\slsf any} sequence space is $e_1=(1,0,...,0,...), 
..., e_k=(0,...,0,k,0,...) $ etc. The ball of radius $\beta$ centered at $\z$ in $\el_{\r}$ will be denoted as 
$B_{\r}^{\infty}(\z,\beta)$. The corresponding topology will be called the standard one or, an $\r$-topology. 

A {\slsf polyradius} is a sequence $\ro =\{ \rho_k\}$ such that all $\rho_k>0$. If $\ro =\{ \rho_k\}\in \el_r$
we call it an $\r$-{\slsf polyradius}. Fix an polyradius $\ro = \{\rho_k\}_{k=1}^{\infty}$. A polydisk of polyradius 
$\ro$ is the set $\Delta^{\infty}(z^0,\ro):=\{\z\in \el_{\r}: |z_k-z_k^0|< \rho_k \text{ for all } k\}$. It is not 
an open subset of $\el_{\r}$ in general. The ``closed'' polydisk is $\bar\Delta^{\infty}(z^0,\ro):=\{ z\in \el_{\r}:
|z_k-z_k^0| \le \rho_k\}$ and it is a closed subset of $\el_{\r}$. For a real positive  $\beta$ and a polyradius 
$\rho$ we set $\beta\rho = \{\beta\rho_n\}_{n=1}^{\infty}$. Note that for the special polyradius $\r$
we have that $\Delta^{\infty}(0,\beta\r) \subset \bar B_{\r}^{\infty} (0,\beta)\subset \Delta^{\infty} (0, (\beta 
+\eps)\r )$ for any $\eps >0$. We also denote $\Delta^{\infty}\left(0, \beta\frac{\adyn}{\n}\right)$ 
as $\Delta^{\infty}\left(0, \frac{\beta}{\n}\right)$. We call the polydisk $\Delta^{\infty}(0,\beta\r)$ 
{\slsf congruent} to $\Delta^{\infty}(0,\r)$, the same for $\Delta^{\infty}(z^0,\beta\r) = z^0 + \Delta^{\infty}
(0,\beta\r)$ with $\Delta^{\infty}(z^0,\r) = z^0 + \Delta^{\infty}(0,\r)$. 

\smallskip By $\n^{-1}$ or by $\frac{\adyn}{\n}$ we denoted the sequence/polyradius 
$\{\frac{1}{n}\}_{n=1}^{\infty}$ and by $\el_{\n^{-1}}$ or by $\el_{\frac{\adyn}{\n}}$ the corresponding
sequence space. We set $\adyn \deff \{1\}_{n=1}^{\infty}$, and therefore $\el^{\infty}$ ( = $\el_{\adyn}$) is 
the sequence space of bounded complex sequences with the norm $\norm{z}_{\infty}\deff \sup_n|z_n|$. We shall also
consider polyradiuses $\r^{\eta}\deff \{r_n^{\eta}\}_{n=1}^{\infty}$ and $\beta\r^{\eta}\deff 
\{\beta r_n^{\eta}\}_{n=1}^{\infty}$ for $\beta , \eta >0$. Notice that these polyradiuses are not in $\el_{\r}$ 
in general.
 
\smallskip 
Now let us turn to the coordinate cross $\C^{\infty}_{\r}\deff \{\{z_n\}_{n=1}^{\infty}\in \el_{\r}: z_n=0 
\text{ for some } n\}$ in $\el_{\r}$.

\begin{lem}
\label{c-cross}
The closure of the cross $\C^{\infty}_{\r}$ in $\r$-topology is equal to 
\begin{equation}
\eqqno(c-cr-cl)
\bar\C^{\infty}_{\r} = \left\{\z\in \el_{\r}:\inf_{\substack{n\in \nn}}
\frac{|z_n|}{r_n}=0\right\},
\end{equation}
and the open set $\el_{\r}\setminus \bar\C^{\infty}_{\r}$ is pathwise
connected.
\end{lem}
\proof Let $\bar\C^{\infty}_{\r}$ be  the set defined by \eqqref(c-cr-cl). Take some $\z\in \el_{\r}\setminus 
\bar\C_{\r}^{\infty}$,  \ie such that $\eps \deff \inf_{\substack{n\in \nn}}\frac{|z_n|}{r_n}>0$. The ball 
$B^{\infty}(\z,\eps)$ obviously doesn't intersect $\C_{\r}^{\infty}$ and this proves that the complement to 
$\bar\C_{\r}^{\infty}$ is 
open. On the other hand take any $\z\in \bar\C_{\r}^{\infty}$. Then for any $\beta >0$ there exists $n$ such that
$\frac{|z_n|}{r_n}<\beta$. That means that the point $\z^0\deff \z - z_ne_n\in \C_{\r}^{\infty}\cap B^{\infty }
(\z ,\beta)$. This proves that $\bar\C^{\infty}_{\r}$, defined by \eqqref(c-cr-cl) is the closure of $\C_{\r}^{\infty}$
in $\el_{\r}$. 

\smallskip To prove the connectivity of $\el_{\r}\setminus \bar\C_{\r}^{\infty}$ take some other 
$\w\in \el_{\r}\setminus \bar\C_{\r}^{\infty}$ with $\delta\deff \inf_{\substack{n\in \nn}}\frac{|w_n|}{r_n}>0$. 
We need to construct a continuous path from $\w$ to $\z$ in $\el_{\r}\setminus \bar\C_{\r}^{\infty}$. 
Consider the following path: in every coordinate go from $w_n$ to $|w_n|$ along the circle of 
radius $|w_n|$ in the clockwise direction, then to $|z_n|$ along the real axis and then along the circle of radius 
$|z_n|$ to $z_n$ in the counterclockwise direction. I.e., $\gamma = \gamma_0 + \gamma_1 + \gamma_2$, where 
$\gamma_0(t)= \{|w_n|e^{(1-t)i\arg w_n}\}$, $\gamma_1(t) = \{|w_n| + t(|z_n|-|w_n|)\}_{n=1}^{\infty}$ and 
$\gamma_2(t) = \{|z_n|e^{ti\arg z_n}\}$. In all cases $t\in [0,1]$ and $\arg$ is supposed to take its values
in $[0, 2\pi )$. Remark that for avery $t\in [0,1]$ one has that 
\[
\frac{|\gamma_0^n(t)|}{r_n} = \frac{|w_n|}{r_n} \ge \delta >0,
\]
which shows that $\gamma_0$ is a path in $\el_{\r}\setminus \bar\C_{\r}^{\infty}$. The same for $\gamma_2$.
As for $\gamma_1$ we see that
\[
\frac{|\gamma_1^n(t)|}{r_n} = \frac{(1-t)|w_n| + t|z_n|}{r_n} \ge \min\{\eps ,\delta\} >0,
\]
which, again proves that $\gamma_1$ is a path in $\el_{\r}\setminus \bar\C_{\r}^{\infty}$. To check the 
continuity of $\gamma_0$ take some $t_0\in [0,1]$ and $t_1\in [0,1]$ and write 
\[
\left||w_n|e^{(1-t_1)i\arg w_n} - |w_n|e^{(1-t_0)i\arg w_n}\right|\le \frac{3}{2}|w_n||t_1-t_0|2\pi .
\]
This implies that $\norm{\gamma_0(t_1) - \gamma_0(t_0)}_{\r} \le 3\pi\norm{w}_{\r}|t_1-t_0|$, 
\ie that $\gamma_0$ is $\r$-continuous. The same work for $\gamma_2$. As for $\gamma_1$ write 
\[
\frac{|\gamma_1^n(t_1)-\gamma_1^n(t_0)|}{r_n} = \left|\frac{(1-t_1)|w_n| + t_1|z_n| - (1-t_0)|w_n| - t_0|z_n|}{r_n}
\right|=\frac{|t_0-t_1|(|w_n| + |z_n|)}{r_n}.
\]
This shows that $\norm{\gamma_1(t_1)-\gamma_1(t_0)}_{\r}\le |t_1-t_0|\left(\norm{z}_{\r} + \norm{w}_{\r}\right)$, 
\ie that $\gamma_1$ is continuous. 

\smallskip\qed

It is easy to see that paths appeared in the proof of this lemma have {\slsf bounded argument}.
By saying this we mean the following:

\begin{defi}
\label{b-arg1}
We say that a continuous path $\gamma (t) = \{\gamma_n(t): t\in [0,1]\}_{n\in \nn}$ in 
$\el_{\r}\setminus \C^{\infty}_{\r}$ is a path with a bounded argument if for every $n\in \nn$
there exists a continuous  choice of $\arg \gamma_n(t)$ such that 
\begin{equation}
\eqqno(b-arg2)
\Pi \deff \Sup_{\substack{t\in [0,1]\\ n\in \nn}}\{|\arg \gamma_n(t)|\}  < + \infty .
\end{equation}
\end{defi}

\begin{lem} 
\label{b-arg2}
Every continuous path in $\el_{\r}\setminus \bar\C^{\infty}_{\r}$ is a path with a bounded argument.
\end{lem}
\proof First let us make a remark about paths in the complement of the cross. Notice that for a point 
$\p=\{p_n\}\in \el_{\r} \setminus \bar\C^{\infty}_{\r}$ 
one has that 
\begin{equation}
\eqqno(dist-cr1)
\dist (\p, \bar\C^{\infty}_{\r}) = \dist (\p, \C^{\infty}_{\r}) = \inf_{n\in \nn} 
\left\{\frac{|p_n|}{r_n}\right\}.
\end{equation}
Indeed, 
\[
\dist (\p, \C^{\infty}_{\r}) = \inf_{n\in \nn}\dist (\p, \{z_n=0\}) = \inf_{n\in \nn}\left\{\inf\{\norm{\p-\z}_{\r}: 
z_n=0\}\right\} = \inf_{n\in \nn}\left\{\frac{|p_n|}{r_n}\right\}.
\]
Therefore for a continuous path $\gamma$ in $\el_{\r}\setminus \bar\C^{\infty}_{\r}$ the distance from $\gamma (t)$
to $\bar\C^{\infty}_{\r}$ is equal to 
\begin{equation}
 \eqqno(dist-cr2)
 \dist (\gamma (t), \bar\C^{\infty}_{\r}) = \inf_{n\in \nn}\left\{\frac{|\gamma_n(t)|}{r_n}\right\},
\end{equation}
and it depends continuously on $t$. In particular there exists $d>0$ such that 
\begin{equation}
 \eqqno(dist-cr4)
\frac{|\gamma_n(t)|}{r_n}\ge d \quad\text{ for all } \quad n.
\end{equation}

\noindent
Now for a continuous path $\gamma :[0,1]\to \el_{\r}\setminus \bar\C^{\infty}_{\r}$ fix the argument of 
every $\gamma_n(0)$ in the interval $[0,2\pi)$ and then extend it continuously in $t\in [0,1]$ for every
$n$. Since $\gamma $ is uniformly continuous we see that for every $\eps >0$ there exists $\delta >0$ such that 
\begin{equation}
\eqqno(b-arg4)
\left|\frac{\gamma_n(t_1)}{r_n}-\frac{\gamma_n(t_2)}{r_n}\right| \le \eps \quad\text{ for all } \quad n,
\end{equation}
provided $|t_1-t_2|\le \delta$. Rewrite this as 
\begin{equation}
\eqqno(b-arg5)
\left|\frac{|\gamma_n(t_1)|e^{i\arg \gamma_n(t_1)}}{r_n}-\frac{|\gamma_n(t_2)|e^{i\arg\gamma_n(t_2)}}{r_n}
\right|\le \eps \quad\text{ for all } \quad n.
\end{equation}
Via \eqqref(dist-cr4) we see that 
\[
\eps \ge \left|\frac{|\gamma_n(t_1)|e^{i\arg \gamma_n(t_1)}}{r_n}-\frac{|\gamma_n(t_2)|e^{i\arg\gamma_n(t_2)}}{r_n}\right|=
\]
\[
= \left|\frac{|\gamma_n(t_1)|}{r_n}\left(e^{i\arg \gamma_n(t_1)}-e^{i\arg \gamma_n(t_2)}\right)
+ \frac{|\gamma_n(t_1)|-|\gamma_n(t_2)|}{r_n}e^{i\arg\gamma_n(t_2)}\right|
\ge 
\]
\[
\ge d\left|e^{i\arg \gamma_n(t_1)}-e^{i\arg \gamma_n(t_2)}\right| - \frac{|\gamma_n(t_1)-\gamma_n(t_2)|}{r_n}. 
\]
Which shows that $|\arg\gamma_n(t_1)-\arg\gamma_n(t_2)|\le \frac{2\eps}{2\pi d} = \frac{1}{\pi d}$ for all $n$ provided 
$|t_1-t_2|\le \delta$. This proves the assertion of the lemma.

\smallskip\qed

\begin{rema} \rm
\label{inter}
a) The (non-open) set $\el_{\r}\setminus \C_{\r}^{\infty}$ is connected as well. To see this take for  a 
given points $\z , \w \in \el_{\r}\setminus \C_{\r}^{\infty}$ the line $l$ joining them. $l\cap (\el_{\r}\setminus 
\C_{\r}^{\infty})$ is either empty (\ie $l$ is contained in $\C_{\r}^{\infty}$) or, is a complement to an at most
countable set $l\cap \C_{\r}^{\infty}$. In our case $l$ is not contained in $\C_{\r}^{\infty}$. Therefore  
$l\cap (\el_{\r}\setminus \C_{\r}^{\infty})$  is connected. This implies the connectivity  (by paths) 
of $\el_{\r}\setminus \C_{\r}^{\infty}$. 

\noindent b) At the same time not every continuous path in $\el_{\r}\setminus \C_{\r}^{\infty}$ is a path with a 
bounded argument. For example $\gamma (t) = \{\frac{1}{n^{2+it}}\}_{n=1}^{\infty}$ is a $\calc^{\infty}$-path 
in $\el_{\n^{-1}}\setminus \C_{\n^{-1}}^{\infty}$ with an unbounded argument.

\noindent c) The intersection $l\cap \C_{\r}^{\infty}$, where $l$ is a complex line such that 
$l\not\subset \C^{\infty}_{\r}$, can be an arbitray countable set in the topology of $l$. Consider, for example, 
the following line $l$ in $\el_{\r}$
\[
z_n = e^{i\theta_n}r_n^2 - r_n\lambda , \quad \lambda\in \cc , \text{ and } \theta_n\in [0,2\pi ) 
\text{ are fixed }.
\]
Here we asume that the sequence $\{r_n\}$ is bounded. Then $l\cap \C_{\r}^{\infty} =
\{r_ne^{i\theta_n}\}_{n=1}^{\infty}$ in the topology of $l$.

\noindent d) Notice also that $\bar\C_{\r}^{\infty}$ contains the line $l_1=\{(\lambda,\{r_n^2\lambda\}_{n=2}^{\infty})
:\lambda\in \cc\}$ entirely (provided $r_n\to 0$), but $l_1\cap C^{\infty}_{\r} =\{0\}$. 
\end{rema}

\begin{lem}
\label{two-crosses}
Let $\rho=\{\rho_n\}_{n=1}^{\infty}$ and $\r= \{r_n\}_{n=1}^{\infty}$ be two polyradiuses such that 
$\frac{\rho_n}{r_n} \to 0$. Then for a point $\z^0 = \{z^0_n\}_{n=1}^{\infty}\in \el_{\r}\setminus 
\bar\C^{\infty}_{\r}$ one has
\begin{equation}
\eqqno(e-c-cr1)
\{\z^0 + \el_{\rho}\}\cap \bar\C^{\infty}_{\r} = \{\z^0 + \el_{\rho}\}\cap\C^{\infty}_{\r}.
\end{equation}
Moreover, the set on the right-hand side of \eqqref(e-c-cr1) is a countable union of complex hyperplanes 
leaving any bounded in the topology of $\el_{\rho}$ subset of the affine plane $\z^0 + \el_{\rho}$.
\end{lem}
\proof For a given $\z^0=\{z^0_n\}\in \el_{\r}\setminus \bar\C^{\infty}_{\r}$ let $\z = \{z_n\}_{n=1}^{\infty}
\in \el_{\rho}$ be such that $\z^0+\z \in \bar\C^{\infty}_{\r}$, \ie 
we have that 
\begin{equation}
\eqqno(e-c-cr2)
\inf_{\substack{n\in \nn}}\frac{|z_n^0+z_n|}{r_n} = 0, \quad\text{ but }\quad
\inf_{\substack{n\in \nn}}\frac{|z_n^0|}{r_n} =: \eps >0.
\end{equation}
Since
\begin{equation}
\eqqno(e-c-cr3)
\beta \deff \sup_{\substack{n\in \nn}}\frac{|z_n|}{\rho_n} < + \infty
\end{equation}
we get
\[
\frac{|z_n^0+z_n|}{r_n} \ge \frac{|z_n^0|}{r_n} - \frac{|z_n|}{\rho_n}\frac{\rho_n}{r_n} 
\ge \eps - \beta\frac{\rho_n}{r_n}.
\]
Due to the assumed $\rho_n/r_n\to 0$ we see that 
\[
\frac{|z_n^0+z_n|}{r_n} \ge \frac{\eps}{2} \quad\text{ for } \quad n\gg 1.
\]
Therefore the first infinum in \eqqref(e-c-cr2) is achieved for a finite $n$, \ie $z^0_n+z_n=0$
for some $n$. This proves that $\z^0 + \z \in \C^{\infty}_{\r}$. The opposite inclusion being trivial 
this proves the identity \eqqref(e-c-cr1).

\smallskip Setting $L_n \deff \{\z=\{z_m\}_{m=1}^{\infty}\in \el_{\rho}: z_n + z^0_n=0\}$ we see that 
\[
\{\z^0 + \el_{\rho}\}\cap \C^{\infty}_{\r} = \bigcup_{n\in \nn}L_n,
\]
is a countable union of hyperplanes. Suppose that an infinite number of $L_n$-s do intersect some 
bounded in the topology of $\el_{\rho}$ subset of the affine hyperplane $\z^0 + \el_{\rho}$. That 
means that $\exists$ $\z^0+ \z^k\in L_{n_{k}}$ with $n_k\to \infty$ such that $\norm{\z^k}_{\rho}
\le C$ for some constant $C$ indpendent of $k$. But $\z^0 + \z^k\in L_{\n_k}$ means that $z^k_{n_k} 
= -z^0_{n_k}$.  Therefore $|z^0_{n_k}|/\rho_{n_k}\le C$.
At the same time $|z^0_{n_k}|/r_{n_k}\ge \eps >0$ because it was assumed that $\z^0\in \el_{\r}\setminus 
\bar\C_{\r}$. This is a contradiction, indeed
\[
\frac{|z^0_{n_k}|}{r_{n_k}} = \frac{|z^0_{n_k}|}{\rho_{n_k}}\frac{\rho_{n_k}}{r_{n_k}}
\le C\frac{\rho_{n_k}}{r_{n_k}} \to 0.
\]

\smallskip\qed

\newprg[VV.zeta]{Riemann zeta function, definition}

For  $\sigma \in \bar\rr$ we had set $H_{\sigma} \deff \{s\in \cc: \re s > \sigma\}$. In particular $H_1 = 
\{z\in \cc:\re s >1\}$, $H\deff H_0=\{s\in \cc: \re s>0\}$, as well as $H_{+\infty} = \emptyset$ and 
$H_{-\infty} = \cc$. For $\b\in \el^{\infty}$, $\z=\{z_n\}\in B^{\infty}_{\n^{-1}} \left(\n^{-1},
1\right)$ and  $s\in H_1$ we had set 

\begin{equation}
\eqqno(zeta-f)
\zeta (\b, \z,s) = \sum_{n=1}^{\infty}b_nz_n^s.
\end{equation}
For the usual Riemann zeta function $\zeta$ one has $\zeta (s) = \zeta (\adyn ,\n^{-1},s)$ 
in our notations.

\begin{lem} 
\label{zeta-def}
Function $\zeta$ is well defined and holomorphic on $\el^{\infty}\times B^{\infty}_{\n^{-1}} (\n^{-1},
1) \times H_1$ as a function of the triple $(\b,\z,s)$.
\end{lem}
\proof Remark that for $z=\{z_n\}_{n=1}^{\infty}\in B^{\infty}_{\n^{-1}} (\n^{-1},1)$ 
one has $\re z_n >0$ for all $n$ and therefore the principal determinations $\Ln z_n = \ln |z_n| + i\arg z_n$ with 
$\arg z_n\in (-\pi /2, \pi/2)$ are well defined. Write now 
\[
\sum_{n=1}^{\infty}|b_nz_n^s| = \sum_{n=1}^{\infty}|b_n||e^{s \cdot \Ln z_n}| = \sum_{n=1}^{\infty}
|b_n|e^{\re s \cdot \log |z_n| - \im s \cdot \arg z_n }. 
\]
Take an exhaustive sequence of compacts in $H_1$, say 
\begin{equation}
\eqqno(comp-k)
K = \{s: 1+1/k\le \re s\le k \text{ and } |\im s |\le k\}.
\end{equation}
Since $|\arg z_n |\le \frac{\pi}{2}$ we see that $e^{- \im s \cdot \arg z_n }\le e^{k\frac{\pi}{2}}$
for all $n\in \nn$ and all $s\in K$. From this fact and from the observation that $|z_n|< \frac{2}{n}$ for $\z\in 
B^{\infty}(\n^{-1},1)$ we obtain 
\begin{equation}
\eqqno(comp-k1)
\sum_{n=1}^{\infty}|b_nz_n^s| \le \norm{b}_{\infty}e^{k\frac{\pi}{2}}\sum_{n=1}^{\infty}
\left(\frac{2}{n}\right)^{1+1/k}
\end{equation}
uniformly on $B\times B^{\infty}(\n^{-1},1)\times K$ for any bounded $B\subset \el^{\infty}$. 
This implies the needed statement. Indeed, the partial sums 
\begin{equation}
\eqqno(part-sum)
\sum_{n=1}^Nb_nz_n^s = \sum_{n=1}^Nb_ne^{s\Ln z_n}
\end{equation}
of \eqqref(zeta-f) are functions of finitely many variables, which are obviously holomorphic, and by 
what was just proved, converge to $\zeta (b, z,s)$ uniformly on $B\times B^{\infty}(\n^{-1},
1)\times K$ for every bounded $B\subset \el^{\infty}$ and compact $K\comp H_1$. Therefore $\zeta (\b , \z ,s)$
is holomorphic, see Proposition \ref{unif} in the Appendix.

\smallskip\qed

Let us give a definition, though obvious, of the notion of analytic continuation in $\el_{\r}$ (in 
any complex Banach space in fact). Let $\gamma :[0,1]\to \el_{\r}$ be a continuous path and let $f_0$
be a holomorphic function (an element) in some  open connected $V_0\ni \gamma (0)$.

\begin{defi}
\label{an-cont}
An analytic continuation of the element $f_0$ along $\gamma$ is a collection of open connected neighborhoods 
$V_t$ of $\gamma (t)$ and $\r$-holomorphic on $V_t$ functions $f_t$ such that for $|t_1 -t_2|$ small
enough one has 
\begin{equation}
\eqqno(an-cont1)
f_{t_1}|_{V_{t_1}\cap V_{t_2}} = f_{t_2}|_{V_{t_1}\cap V_{t_2}},
\end{equation}
\end{defi}

\begin{lem}
\label{l2-cont}
Function $\zeta$ admits an analytic continuation along any continuous path in 
\begin{equation}
\el^{\infty}\times \left(\el_{\n^{-1}}\setminus \bar\C^{\infty}_{\n^{-1}}\right)\times H_1.
\end{equation}
\end{lem}
\proof Along the proof of the lemma we shall assume that $b$ is fixed in some bounded open $B\subset \el^{\infty}$ 
and $s$ is fixed in $K\deff \{1+1/k\le\re s\le k, |\im s|\le k\}$. Take a path  $\gamma = \{\gamma (t): t\in [0,1]\}$ 
from the distinguished point $\n^{-1}$ to some $\z = \{z_n\}$ in $\el_{\n^{-1}}\setminus \bar \C^{\infty}_{\n^{-1}}$ 
with the argument bounded by $\Pi$ as in Definition \ref{b-arg1}. Due to \eqqref(dist-cr2) the distance from 
$\gamma (t)$ to $\bar\C^{\infty}_{\n^{-1}}$ is separated from zero, say by $\eps >0$. Therefore 
\begin{equation}
\eqqno(dist-cr3)
|\gamma_n(t)|\ge \frac{\eps}{n} \quad\text{ for all } \quad n\in \nn.
\end{equation}

\smallskip\noindent{\slsf Step 1.} {\it Continuation along $\gamma $.} As open sets take $V_t\deff B\times 
B^{\infty}_{\n^{-1}} (\gamma (t), \eps)\times H_1$.
%where polyradius $\ro_t$ is defined as
%\[
%\ro_t = |\gamma (t)| \deff \left\{|\gamma_n (t)|\right\}_{n\in \nn }.
%\]
One obvioulsly has that $V_t\subset \el_{\n^{-1}}\setminus \bar \C^{\infty}_{\n^{-1}}$ for every $t$ due to
\eqqref(dist-cr3), as well as $V_{t_1}\cap V_{t_2}\not=\emptyset$ for $t_1,t_2\in [0,1]$ close enough. 
Function $\zeta_t(\b,\w,s)$ for $(\b, \w,s)\in V_t$ define naturally as follows. Notice first that every
$\w=\{w_n\}\in B^{\infty}(\gamma (t), \eps)$  can be uniquely written in the form 
\begin{equation}
\eqqno(zeta-t1)
w_n = \gamma_n(t)\left[1+r_ne^{i\theta_n}\right] \quad\text{ for } n=1,...
\end{equation}
with some $\theta_n\in [0,2\pi)$ and some $0\le r_n<1$. The latter is because one should have $|nr_n\gamma_n(t)| <\eps$.
We set
\begin{equation}
\eqqno(zeta-t)
\zeta_t(\b,\w,s) = \sum_{n=1}^{\infty} b_nw_n^s = \sum_{n=1}^{\infty}b_ne^{s\ln w_n},
\end{equation}
were the arguments are defiend as $\arg w_n = \arg \gamma_n(t) + \arg \left[1+r_ne^{i\theta_n}\right]$.
The first summand was predefined by the boundedness of argument condition imposed on $\gamma$ and the second
varies in $]-\pi/2,\pi/2[$ and is determined uniquely by $r_n$ and $\theta_n$. Remark that series \eqqref(zeta-t) 
converge normally on $B\times B^{\infty}_{\n^{-1}}(\gamma (t), \eps)\times K$. Indeed, for $\b\in B$, $\w\in
B^{\infty}_{\n^{-1}}(\gamma (t), \eps))$ and $s\in K\deff \{1+1/k\le\re s\le k, |\im s|\le k\}$ our series  
can be estimated as follows
\[
\sum_{n=1}^{\infty}\left|b_n\gamma_n(t)^s\left[1+r_ne^{i\theta_n}\right]^s\right| \le 
\norm{b}_{\infty}e^{(\Pi +\pi/2)k}\sum_{n=1}^{\infty}\left|\gamma_n(t) \left[1+r_ne^{i\theta_n}\right]\right|^{\re s}
\le
\]
\begin{equation}
\eqqno(est-l+)
\le  \norm{b}_{\infty} e^{(\Pi + \pi)k}\sum_{n=1}^{\infty}\left|\gamma_n(t)\right|^{\re s} \le \norm{b}_{\infty}e^{(\Pi +\pi )k}
\norm{\gamma (t)}_{\n^{-1}}^{\re s}\sum_{n=1}^{\infty}\frac{1}{n^{\re s}}.
\end{equation}
Here we used the idenity $|a^s| = |a|^{\re s}e^{-\im  s \arg a}$ and the bound \eqqref(b-arg2) imposed on 
the argument of $\gamma$. Remark that the estimate \eqqref(est-l+) is uniform along $B\times B^{\infty}
(\gamma (t), \eps)\times K$ because the $\el_{\n^{-1}}$-norm of $\gamma (t)$ is bounded. This proves that 
$\zeta_t$ is holomorphic in  $V_t$. 

\smallskip\noindent{\slsf Step 2.} The fact that $\zeta_{t_1}|_{V_{t_1}\cap V_{t_2}} = 
\zeta_{t_2}|_{V_{t_1}\cap V_{t_2}}$ is obvious since both are defined by the same series \eqqref(zeta-t) 
and the arguments of $w_n$ do not depend on $t_1\sim t_2$. The latter is again due to the 
possibility to choose each $\arg \gamma_n (t)$ continuously and uniformly bounded on $n$. Lemma is proved.

\smallskip\qed

\begin{rema} \rm
We know that every $\z\in \el_{\n^{-1}}$ can be joined with $\n^{-1}$ by a path with bounded argument. 
Moreover, consider the following path from $\frac{1}{\n}$ to $\z$: in every coordinate 
go from $\frac{1}{n}$ first to $|z_n|$ along the real axis and then along the circle of radius $|z_n|$ to 
$z_n$ in the counterclockwise direction. I.e., $\gamma = \gamma_1 + \gamma_2$, where $\gamma_1(t) = 
\{\frac{1}{n} + t(|z_n|-\frac{1}{n})\}_{n=1}^{\infty}$ and $\gamma_2(t) = \{|z_n|e^{ti\arg z_n}\}$, in both 
cases $t\in [0,1]$. $\arg$ can be supposed to vary in $[0, 2\pi )$ only. 
\end{rema}

\newprg[VV.taylor]{Taylor expansion of function zeta}

For $(\b, \z, s) \in \el^{\infty}\times B^{\infty}_{\n^{-1}}(\n^{-1},\adyn)\times H_1$ 
we have 
\[
\zeta (\b, \z, s) = \sum_{n=1}^{\infty} b_n z_n^s,
\]
and therefore we can easily compute all partial derivatives of $\zeta$ with respect to $\z$ to get
\[
\frac{\d^m\zeta}{\d z_n^m}(\b,\z,s) = s(s-1)...(s-m+1)b_nz_n^{s-m} \quad \text{ for all } n,m\in \nn.
\]
Corresponding $m$-homogeneous polynomial in the Taylor expansion of $\zeta$ at $\z$ is 
\begin{equation}
\eqqno(pol-m-zt)
P_m(\b, \z, \w, s) = \frac{s(s-1)...(s-m+1)}{m!}\sum_{n=1}^{\infty}b_nz_n^{s-m}w_n^m.
\end{equation}
This gives the Taylor expansion of zeta function at $(\b,\z,s)\in \el^{\infty}\times 
B^{\infty}_{\n^{-1}} (\n^{-1},\adyn)\times H_1$ with respect to the variable $\z$:
\begin{equation}
\eqqno(tay-z)
\zeta (\b, \z+\w, s) - \zeta (\b, \z,s) = \sum_{m=1}^{\infty}\frac{s(s-1)...(s-m+1)}{m!}
\sum_{n=1}^{\infty}b_nz_n^{s-m}w_n^m.
\end{equation}
%Writing this formulas at $z=\n^{-1}$ and get
%\begin{equation}
%\eqqno(tay-n)
%\zeta \left(\adyn , \n^{-1} + w, s\right) - \zeta \left(b,\n^{-1}, s\right) = \sum_{m=1}^{\infty}
%\frac{s(s-1)...(s-m+1)}{m!}\sum_{n=1}^{\infty}\frac{b_n}{n^{s-m}}w_n^m.
%\end{equation}
It occurs that the right hand side of \eqqref(tay-z) makes sense and converge for $\z =\{z_n\}\in 
\el_{\n^{-1}} \setminus \bar C^{\infty}_{\n^{-1}}$ with the silent assumption that arguments of 
$z_n$ are chosen to be bounded, \ie 
\begin{equation}
\eqqno(b-arg3)
\Pi \deff \sup \{|\arg z_n|\} < \infty.
\end{equation}
For example $\z$ is an endpoint of a path with bounded argument.
\begin{lem}
\label{taylor-z}
For the series on the right hand side of \eqqref(tay-z) the following holds true.

\smallskip\noindent \sli For every $(\b,\z,s)\in \el^{\infty}\times \left(\el_{\n^{-1}}\setminus 
\bar\C_{\n^{-1}}\right)\times H_1$ all polynomials $P_m\left(b, \z + w, s\right)$ are 
$\n^{-1}$-bounded 

\quad and the  series \eqqref(tay-z) uniformly converge on $\{\text{ bounded sets in } \el^{\infty}\}
\times \bar B^{\infty}_{\n^{-1}}(0, \beta )\times$

\quad  $\{\text{compacts of } H_1\}$ for every $0<\beta<\norm{\z}_{\n^{-1}}$.

\smallskip\noindent \slii For every $0 \le \eta < + \infty $ and every $0<\beta<\norm{\z}_{\n^{-1}}$ the power 
series  expansion  \eqqref(tay-z) 

\quad  converge uniformly on $\{\text{ bounded sets in } \el^{\infty}\}\times \bar\Delta^{\infty}\left(0, 
\frac{\beta}{\n^{1+\eta}}\right) \times \{\text{compacts of } H_{1-\eta}\}$.

\smallskip\noindent \sliii For all $0 < \eta < + \infty $ and $0 < \beta < \norm{\z}_{\n^{-1}}$ the power series 
expansion \eqqref(tay-z) 
converge 

\quad uniformly on $\{\text{ bounded sets in } \el^{\infty}\}\times \bar\Delta^{\infty}\left(0, 
\frac{\beta e^{-\eta\n}}{\n}\right)
\times \{\text{compacts of } \cc\}$.
\end{lem}

\proof Here, as usual, $H_{1-\eta} \deff \{\re s > 1 - \eta\}$ and  $\frac{\beta}{\n}$ 
is the polyradius $\{\frac{\beta}{n}\}_{n=1}^{\infty}$. By $\frac{\beta}{\n^{1-\eta}}$ we denote the polyraidus 
$\{\frac{\beta}{n^{1-\eta}}\}_{n=1}^{\infty}$ and by $\frac{\beta e^{-\eta\n}}{\n}$ we denote the polyradius 
$\{\frac{\beta e^{-\eta n}}{n}\}_{n=1}^{\infty}$.

\smallskip\noindent (\sli  As for the $\n^{-1}$-boundedness of $P_m$ fix $(\b,\z)\in \el^{\infty}\times 
\left(\el_{\n^{-1}}\setminus \bar\C_{\n^{-1}}\right)$, $\w\in \bar B^{\infty}_{\n^{-1}}(0, \beta )$ and write,
taking into account that $|z_n|\le \norm{z}_{\n^{-1}}/n$ and $|w_n|< \beta /n$ for all $n$, that
\begin{equation}
\eqqno(tay-1)
\left|\sum_{n=1}^{\infty}b_nz_n^{s-m}w_n^m\right| \le \sum_{n=1}^{\infty}\frac{|b_n|
\norm{z}_{\n^{-1}}^{\re s-m}e^{k\Pi}}{n^{\re s-m}}|w_n^m| \le \norm{b}_{\infty}e^{k\Pi}\sum_{n=1}^{\infty}
\frac{\beta^m\norm{z}_{\n^{-1}}^{\re s-m}}{n^{\re s}} \le  
\end{equation}
\[
\le \norm{b}_{\infty}e^{k\Pi}
\norm{z}_{\n^{-1}}^{\re s -m}\beta^m \sum_{n=1}^{\infty}\frac{1}{n^{\re s}} = 
e^{k\Pi}\norm{b}_{\infty}\norm{z}_{\n^{-1}}^{\re s -m}\beta^m \zeta (\re s).
\]
Here $s$ is taken in the compact $\{1 + \frac{1}{k}\le \re s \le k, -k\le \im s \le k\}$.
When $\beta =1$ this gives the  $\n^{-1}$-boundedness of $P_m\left(b, \z, w, s\right)$. As for the
convergence observe first the following Taylor expansion of the holomorphic function $\lambda^s$ at $\adyn$
\begin{equation}
\eqqno(tay-s1)
\lambda^s = 1 + \sum_{m=1}^{\infty}\frac{s(s-1)...(s-m+1)}{m!}(\lambda-1)^m.
\end{equation}
Since $\lambda^s$ is holomorphic in $H=\{\re \lambda >0\}$ the radius of convergence of series in 
\eqqref(tay-s1) is $\adyn$. In particular series 
\begin{equation}
\eqqno(tay-s2)
\sum_{m=1}^{\infty}\frac{|s(s-1)...(s-m+1)|}{m!}t^m
\end{equation}
converge absolutely for any $|t|<1$. For $\w\in \bar\Delta^{\infty}(0,\frac{\beta}{\n})$ we get from \eqqref(tay-1)
that in \eqqref(tay-z)  one has
\begin{equation}
\eqqno(tay-2)
\sum_{m=1}^{\infty}\left|P_m(\b,\z, \w, s)\right|\le \zeta (\re s)\norm{b}_{\infty}
\sum_{m=1}^{\infty}\frac{|s(s-1)...(s-m+1)|}{m!}\beta^m{\norm{z}_{\n^{-1}}}^{\re s -m} =
\end{equation}
\[
= \zeta (\re s)\norm{b}_{\infty}\norm{z}_{\n^{-1}}^{\re s}
\sum_{m=1}^{\infty}\frac{|s(s-1)...(s-m+1)|}{m!}\left(\frac{\beta}{\norm{z}_{\n^{-1}}}\right)^m.
\]
Estimate \eqqref(tay-2) implies, provided $\beta < \norm{\z}_{\n^{-1}}$, the normal convergence of 
\eqqref(tay-z) on  

\noindent $\{\text{ bounded sets in }\el^{\infty}\}\times \bar B^{\infty}_{\n^{-1}}(0, \beta )\times
\{\text{ compacts in } H_1\}$ as claimed.

\smallskip\noindent (\slii Every $\w =\{w_n\}_{n=1}^{\infty}\in \Delta^{\infty}\left(0, \frac{\beta }{\n^{1+\eta}}
\right)$ can be written as $w_n = \frac{\beta r_n}{n^{1+\eta}}e^{i\theta_n}$, where $0\le r_n\le 1$ 
and $\theta_n\in [0,2\pi)$. Now \eqqref(tay-z) writes as 
\[
\zeta \left(\b,\z + \w, s\right) - \zeta (\b, \z, s) = \sum_{m=1}^{\infty}
\frac{s(s-1)...(s-m+1)}{m!}\sum_{n=1}^{\infty}b_nz_n^{s-m}w_n^m = 
\]
\begin{equation}
\eqqno(tay-r)
= \sum_{m=1}^{\infty}\frac{s(s-1)...(s-m+1)}{m!}\beta^m\sum_{n=1}^{\infty}
\frac{b_nr_n^me^{im\theta_n}}{n^{(1+\eta )m}}z_n^{s-m}.
\end{equation}
Since $\z\in \el_{\n^{-1}}\setminus \bar\C^{\infty}_{\n^{-1}}$ we see that for $\re s > 1 - \eta$ and 
$\w\in \bar\Delta^{\infty}\left(0, \frac{\beta}{\n^{1+\eta}}\right)$ \eqqref(tay-r) can be estimated as
follows
\[
\sum_{m=1}^{\infty}\frac{|s(s-1)...(s-m+1)|}{m!}\beta^m\sum_{n=1}^{\infty}
\frac{|b_n|\norm{\z}_{\n^{-1}}^{\re s-m}e^{k\Pi}}{n^{\re s+m\eta}} =
\]
\[
= e^{k\Pi}\norm{\z}_{\n^{-1}}^{\re s}\underbrace{\sum_{m=1}^{\infty}
\frac{|s(s-1)...(s-m+1)|}{m!}\frac{\beta^m}{\norm{\z}_{\n^{-1}}^m}}_{I}\cdot 
\underbrace{\sum_{n=1}^{\infty}\frac{|b_n|}{n^{\re s+m\eta}}}_{II}.
\]
Since $\re s + m\eta >1$ for $m\ge 1$ the term {\sl (II)} can be majorated by 
$\norm{b}_{\infty}\zeta (\re s + m\eta)$. Term {\sl (I)} can be estimated by the remark about power series 
\eqqref(tay-s2) since it was supposed that $\beta <\norm{\z}_{\n^{-1}}$. Therefore sums {\sl (I)} 
and {\sl (II)} converge on compacts in $H_{1-\eta}$ when $0<\beta < \norm{\z}_{\n^{-1}}$.

\smallskip\noindent (\sliii   Every $\w$ in the polydisk $\bar\Delta^{\infty}\left(0, \frac{\beta e^{-\eta\n}}{\n}\right)$
writes as $w_n = \frac{\beta r_n e^{-\eta n}}{n}e^{i\theta_n}$ with $0\le r_n\le 1$, $\theta_n\in [0,2\pi)$. Therefore 
\eqqref(tay-z) writes as 
\begin{equation}
\eqqno(tay-r1)
\zeta \left(\b, \z + \w, s\right) - \zeta (\b , \z, s) = \sum_{m=1}^{\infty}\frac{s(s-1)...(s-m+1)}{m!}\beta^m
\sum_{n=1}^{\infty}\frac{b_nr_n^me^{im\theta_n}}{n^{m}}e^{-\eta nm}z_n^{s-m}.
\end{equation}
Since $|z_n|\le \norm{\z}_{\n^{-1}}/n$ the right-hand side of \eqqref(tay-r) can be estimated by 
\begin{equation}
\eqqno(tay-r2)
\norm{\z}_{\n^{-1}}^{\re s}\sum_{m=1}^{\infty}\frac{|s(s-1)...(s-m+1)|}{m!}
\left(\frac{\beta}{\norm{\z}_{\n^{-1}}}\right)^m\sum_{n=1}^{\infty}\frac{e^{-\eta nm}}{n^{\re s}}.
\end{equation}
And the second term in \eqqref(tay-r2) converge normally on compacts in $\cc$ since 
\[
\sum_{n=1}^{\infty}\frac{1}{n^{\re s}}e^{-\eta nm} \le \sum_{n=1}^{\infty}\frac{1}{n^{\re s}}e^{-\eta n}
\]
for all $m\ge 1$. Convergence of \eqqref(tay-r2) follows now from \eqqref(tay-s2).

\smallskip\qed

\newprg[VV.proof]{Proof of Theorem \ref{m-thick}}

Now we turn to the proof of Theorem \ref{m-thick} from the Introduction. First let us state the item (\sliii
of Lemma \ref{taylor-z} in a suitable form. Notice once more that writing that $(\b , \z)\in \calm$ we mean
not only that $\zeta (\b , \z, \cdot)$ extends meromorphically to $\cc$ but also that $\z \in \el_{\n^{-1}}
\setminus \bar\C^{\infty}_{\n^{-1}}$.

\begin{corol}
\label{taylor-c}
For every  $(\b , \z) \in \calm $ the power series expansion 
\begin{equation}
\eqqno(tay-z1)
\zeta (\b , \z + \w,s) - \zeta (\b , \z , s) = \sum_{m=1}^{\infty}\frac{s(s-1)...(s-m+1)}{m!}
\sum_{n=1}^{\infty}b_nz_n^{s-m}w_n^m
\end{equation}
converge uniformly on $B^{\infty}_{e^{-\n}}(0,\beta)\times \{\text{compacts of } \cc \}$ for 
any $0<\beta <\norm{\z}_{\n^{-1}}$. Moreover, this expansion provides the analytic continuation 
of $\zeta (\b, \cdot , \cdot)$ to $\{\z + B^{\infty}_{e^{-\n}}(0,\beta)\}\times \cc $.
\end{corol}
\proof Indeed, notice that for $0<\beta <\norm{\z}_{\n^{-1}}$ we have that 
\begin{equation}
\eqqno(incl-1)
B^{\infty}_{e^{-\n}}(0,\beta) \subset \bar\Delta^{\infty}\left(0,\frac{\beta e^{-\n}}{\n}\right)
\end{equation}
in the topology of $\el_{\e^{-\n}}$. Therefore item (\sliii of Lemma \ref{taylor-z} gives us 
the uniform convergence of \eqqref(tay-z1) on $B^{\infty}_{e^{-\n}}(0,\beta)\times \{\text{compacts of } \cc \}$.
This results to the analytic continuation of the power expansion in question to $B^{\infty}_{e^{-\n}}(0,\beta)
\times \cc $. Since $\zeta (\b, \z , \cdot)$ was supposed to extend meromorphically to $\cc$ we see that
for every $\z + \w\in \{\z + B^{\infty}_{e^{-\n}}(0,\beta)\}$ function $\zeta (\b , \z + \w , \cdot)$ extends
meromorphically to $\cc$. 

\smallskip Now recall that $\zeta (\b, \cdot ,  \cdot)$ is holomorphic on $\left(\el_{\n^{-1}}\setminus 
\bar C^{\infty}_{\n^{-1}}\right)\times H_1 \supset \{\z + B^{\infty}_{e^{-\n}}(0,\norm{\z}_{\n^{-1}})\}\times H_1$.
Remark that (a local) restriction of a $\n^{-1}$-holomorphic function to $\el_{e^{-\n}}$ is $e^{-\n}$-holomorphic.
This follows from the obvious G\^ateaux diffferentiability and continuity of such restriction. 
Therefore we can apply Lemma \ref{siciak} to $\zeta (\b, \z+\w,s)-\zeta (\b, \z,s)$ and conclude that it is
holomorphic on $\{\z + B^{\infty}_{e^{-\n}}(0,\norm{\z}_{\n^{-1}})\}\times \cc$ as a function of a couple $(\z, s)$.
This gives us the meromophicity of $\zeta (\b, \z+\w,s)$ on  $\{\z + B^{\infty}_{e^{-\n}}(0,\beta)\}\times \cc $
as stated.

\smallskip\qed 

\smallskip We continue the proof of Theorem \ref{m-thick}. Let a point $(\b ,\z_0 )\in \calm$ be given, \ie 
$\zeta (\b , \z_0 , s)$ extends in $s$ to a meromorphic function on $\cc$. Then by Corollary \ref{taylor-c} 
the difference $\zeta (\b, \z, s) - \zeta (\b,\z_0,s)$ extends holomorphically to $\left\{\z_0 + 
B^{\infty}_{e^{-\n}}(0,\beta)\}\right\}\times \cc$ for any $0<\beta <\norm{\z}_{\n^{-1}}$.
Take some $\z_1\in \left\{\z_0 + \el_{\e^{-\n}}\right\} \setminus \bar\C^{\infty}_{\n^{-1}}$ and take some 
path $\{\z (t): t\in [0,1]\}$ in this domain from $\z_0$ to $\z_1$ with bounded argument. Notice that by 
Lemma \ref{l2-cont} function $\zeta (\b,\cdot,\cdot)$ admits an analytic continuation along this path
with the restriction that $s\in H_1$. More precisely, we have meromorphic functions $\zeta_t(\b , \cdot , \cdot)$
in $\left\{\z_t + B^{\infty}_{\e^{-\n}}(0,\beta)\right\}\times H_1$ such that $\zeta_{t_1}= \zeta_{t_2}$
on the  intersection $\left\{\z_{t_1} + B^{\infty}_{\e^{-\n}}(0,\beta)\right\}\cap 
\left\{\z_{t_2} + B^{\infty}_{\e^{-\n}}(0,\beta)\right\}\times H_1$ for $|t_1-t_2|$ small enough and  
$\zeta_0=\zeta$ on $\left\{\z_0 + B^{\infty}_{\e^{-\n}}(0,\beta)\right\}\times H_1$. But $\zeta_0$ extends
to  $\left\{\z_0 + B^{\infty}_{\e^{-\n}}(0,\beta)\right\}\times \cc$. At the same time for $t$ close to zero
$\left\{\z_t + B^{\infty}_{\e^{-\n}}(0,\beta)\right\}\cap \left\{\z_0 + B^{\infty}_{\e^{-\n}}(0,\beta)\right\}$ 
is plurithick in $\left\{\z_t + B^{\infty}_{\e^{-\n}}(0,\beta)\right\}$. Therefore by Lemma \ref{siciak} 
$\zeta_t$ extends to $\left\{\z_t + B^{\infty}_{\e^{-\n}}(0,\beta)\right\}\times \cc$. The rest is obvious and
our Theorem is proved.

\smallskip\qed

\newprg[VV.dir]{Dirichlet series}

Let us interpret the item (\sliii of Lemma \ref{taylor-z} in terms of the general Dirichlet series
\begin{equation}
\eqqno(dir-gen)
\cald (\a , \lambda , s) \deff \sum_{n=1}^{\infty}a_ne^{-\lambda_ns} = \sum_{n=1}^{\infty}\frac{a_n}{e^{\lambda_ns}},
\end{equation}
where  $\a=\{a_n\}_{n=1}^{\infty}\in \el^{\infty}$ and $\lambda = \{\lambda_n\}_{n=1}^{\infty}$ are  complex. We 
require that $\sup |ne^{-\lambda_n}| <\infty$, \ie that $\{e^{-\lambda_n}\}\in \el_{\n^{-1}}$.
If, for example, $a_n=1$ then we deal with a Dirichlet series tout court
\begin{equation}
\eqqno(dir-gen1)
\cald (\lambda , s) \deff \sum_{n=1}^{\infty}e^{-\lambda_ns} = \sum_{n=1}^{\infty}\frac{1}{e^{\lambda_ns}}.
\end{equation}
Now (\sliii means that for $z_n= e^{-\lambda_n}$ one has that series \eqqref(dir-gen1) meromorphically extend 
to the whole of $\cc$ provided that for some $0<\eta <\infty$ and $0< \beta <1$ one has 
\begin{equation}
\eqqno(dirich1)
\left|e^{-\lambda_n}-\frac{1}{n}\right|<\frac{\beta e^{-\eta n}}{n} \quad\text{ for all }\quad n.
\end{equation}
Even better:
\begin{corol}
\label{dirichlet}
Suppose that \eqqref(dir-gen)  extends meromorphically to $\cc$ for a given $\{\lambda_n\}_{n=1}^{\infty}$. 
Then 
\begin{equation}
\eqqno(dir-gen2)
\cald (\a , \mu ,s) \deff \sum_{n=1}^{\infty}a_ne^{-\mu_ns}
\end{equation}
will extend to $\cc$ for all $\{\mu_n\}_{n=1}^{\infty}$ such that for some $0<\beta ,\eta <\infty$ one has
\begin{equation}
\left|e^{-\lambda_n}-e^{-\mu_n}\right|<\frac{\beta e^{-\eta n}}{n} \quad\text{ for all }\quad n.
\end{equation}
\end{corol}
I.e., the set of $\lambda$ such that $\cald (\a , \lambda , s)$ is meromorphic on $\cc$ is open (and non-empty) in an 
appropriate sequence space. Remark, in addition, that according to \eqqref(abs-c) the abscissa of convergence of 
\eqqref(dir-gen1) is $\adyn$ for the classical case when $\lambda_n = \ln n$.

\begin{rema} \rm 
Recall that for a general Dirichlet series \eqqref(dir-gen1) the abscissa of convergence is given by 
\begin{equation}
\eqqno(abs-c)
\text{ either } \sigma_c=\limsup_{\substack{n\to \infty}}\frac{\ln|a_1+...+a_n|}{\lambda_n} \quad\text{ or, } \quad
\sigma_c=\limsup_{\substack{n\to \infty}}\frac{\ln|a_{n+1}+a_{n+2} + ...|}{\lambda_n}
\end{equation}
depending on either $\sum_{n=1}^{\infty}a_n$ diverge or, converge.
\end{rema}

\begin{rema} \rm
\label{poles-move}
1. Remark that the item (ii) of Lemma \ref{taylor-z} is uniform on $z\in \Delta^{\infty}\left(\frac{\adyn}{\n}, 
\frac{\beta}{\n^2}\right)$. Therefore the function $\zeta_0(z,s)$ is continuous as a function of all variables on 
$\Delta^{\infty}\left(\frac{\adyn}{\n}, \frac{\beta}{\n^2}\right)\times H$ in the standard sense.

\noindent
2. When $z$ moving in $\Delta^{\infty}\left(\frac{\adyn}{\n}, \frac{\beta}{\n^2}\right)$ the pole $s=1$ of $\zeta (z,s)$ 
doesn't move. But it moves in general. Let $\eta \sim 1$ then 
\[
\zeta \left(\frac{1}{\n^{\eta}}, s\right) = \zeta (\n, \eta s) = \frac{1}{\eta s -1} + \zeta_0\left(\frac{\adyn}{\n}, s\right).
\]
I.e., the pole of $\zeta \left(\frac{1}{\n^{\eta}},s\right)$ is $\frac{1}{\eta}$ with residue $\frac{1}{\eta}$.
\end{rema}

\newsect[LL]{Zeta function in the form of Lerch-Lipschitz}

\newprg[LL.ll]{Zeta function in the form of Lerch-Lipschitz}

By $\bar\el^{\infty}_{\im^+} = \{\{a_n\}\in \el^{\infty}, \im a_n \ge 0\}$ let us denote the closure 
of $\el^{\infty}_{\im^+}$, and by $\bar \el^{\infty}_{\re^+}=\{\xi_n\}\in \el^{\infty}: \re \xi_n \ge 0\}$, 
the closure of $\el^{\infty}_{\re^+}$. Having made the substitution in zeta function $\zeta (\b, \z,s)$ as in 
\eqqref(chg-var) we obtained  the following Lerch-Lipschitz form of zeta function
\begin{equation}
\eqqno(ll-f1)
\zeta^{LL}(\a, \xi ,s) = \sum_{n=1}^{\infty}\frac{e^{2\pi in a_n}}{(n + \xi_n)^s}.  
\end{equation}
One immediately sees  that this series converge absolutely for $(\a,\xi)\in \bar\el^{\infty}_{\im^+}\times 
\bar\el^{\infty}_{\re^+}$ and $s\in H_1$. Therefore $\zeta^{LL}$ is well defined, holomorphic for $(\a,\xi , s)\in 
\el^{\infty}_{\im^+}\times \el^{\infty}_{\re^+}\times H_1$ and continuous up to $\bar\el^{\infty}_{\im^+}\times 
\bar \el^{\infty}_{\re^+}\times H_1$. Indeed, we see that $\b=\{b_n\}\in \el^{\infty}$ if $\a\in \bar\el^{\infty}_{\im^+}$.
As for $\xi = \{\xi_n\} \in \bar\el^{\infty}_{\re^+}$ we see that 
\[
|z_n|\le \frac{1}{n + \re \xi_n} \le  \frac{1}{n},
\]
\ie $\z=\{z_n\}\in \el_{\frac{1}{\n}}$. Therefore $\zeta^{LL}$ is obtained 
as a composition of $\zeta$ with mapping \eqqref(chg-var) on the whole domain in question. 
Notice that the ``cross''
\begin{equation}
\eqqno(cross-2)
\C^{\infty}_-\deff \{\z\in \el^{\infty}: z_n = -n \text{ for some } n\} = \bigcup_{n=1}^{\infty}\{z_n=-n\}, 
\end{equation}
is a closed subset of $\el^{\infty}$.
Set $\el^{\infty}_{\zz^-}\deff \deff \el^{\infty}\setminus \C^{\infty}_-$. Let us call a 
set $R$  in $\el^{\infty}_{\zz^-}$ bounded if there exist constants $E,\eps >0$ such that for every $\xi = \{\xi_n\}\in R$ 
one has 
\begin{equation}
\eqqno(b-set1)
\frac{n}{E} \le |n+\xi_n|\le \frac{n}{\eps} \quad \forall n.
\end{equation}
This relation is equivalent to 
\begin{equation}
\eqqno(b-set2)
\frac{\eps}{n} \le \frac{1}{|n+\xi_n|}\le \frac{E}{n} \quad \forall n,
\end{equation}
\ie to the fact that \eqqref(chg-var) maps $\xi$ to $z$ in a bounded subset of 
$\el_{\frac{\adyn}{\n}}\setminus \bar\C^{\infty}_{\frac{\adyn}{\n}}$. As the result from Lemma \ref{l2-cont}
we obtain the following:
\begin{corol}
\label{l2-cont1}
Function $\zeta^{LL}$ admits a multivalued analytic continuation to any bounded open subset of 
\begin{equation}
\el^{\infty}_{\im^+}\times \el^{\infty}_{\zz^-}\times H_1.
\end{equation}
The ramification takes place in variable $\z$ at points of $\zz^-$, more precisely around hypersurfaces of
the kind $\el^{\infty}\times \{-n\}\times \cc$, here $n\in \nn$. With respect to the parameter $\a\in 
\el^{\infty}_{\im^+}$ function $\zeta^{LL}$ is continuous up to $\bar\el^{\infty}_{\im^+}$.
\end{corol}

From now on we rename variable $\xi$  back to $\z$. It is worth to notice that $\zeta^{LL}(\adyn ,\z ,\cdot)$  
can be continued to $H_0$ in variable $s$ whatever $\z\in \el^{\infty}_{\zz^-}$ is. Moreover, the 
following is true.

\begin{thm}
\label{cont-to-h0}
\sli Let $\a^0 =\{a^0_n\}_{n=1}^{\infty}\in \bar\el^{\infty}_{\im^+}$ be such that for some 
$\z^0=\{z^0_n\}_{n=1}^{\infty}\in \el^{\infty}_{\zz^-}$
function $\zeta (\a^0,\z^0,\cdot)$ can be meromorphically continued in variable $s$ to $H_0$. 
Then function $\zeta (\a^0 , \cdot ,\cdot)$ can be meromorphically 
continued to $\el^{\infty}_{\zz^-}\times H_0$ as a (multivalued) function of the couple $(\z ,s)$.
\end{thm}
\proof Take any $\z^1$ in a convex neighborhood of $\z^0$  which doesn't intersects $\C^{\infty}_-$.
Notice that $\frac{\d \zeta}{\d z_n} = -\frac{se^{2\pi i na^0_n}}{(n+z_n)^{s+1}}$. Therefore for 
the straight intervall $\{\z (t)=(1-t)\z^0+t\z^1:t\in [0,1]\}$ from $\z^0$ to $\z^1$ we have 
\begin{equation}
\eqqno(d1-zeta)
\zeta (\a^0, \z^1,s) - \zeta (\a^0, \z^0, s) = \int\limits_0^1\frac{d\zeta (a^0, \z (t), s)}{dt}dt =
-s\int\limits_0^1\sum_{n=1}^{\infty}\frac{e^{2\pi i n a^0_n} (z^1_n-z^0_n)}{(n+z_n(t))^{s+1}}dt.
\end{equation}
Vector $\{e^{2\pi i n a^0_n}(z^1_n-z^0_n)\}_{n=1}^{\infty}$ belongs to $\el^{\infty}$. 
And since $\eps/n\le 1/|n+z_n(t)|\le E/n$ as in \eqqref(b-set2) we see that the sum under 
the integral on the right hand side of \eqqref(d1-zeta) uniformly converges on compacts of 
$\{\re (s +1) > 1\}$, \ie on $H_0$. Due to the assumption of our lemma $\zeta (\a^0, \z_0, s)$ 
can be meromorphically continued to $H_0$. Therefore the same holds for $\zeta (\a^0, \z_1, s)$ 
whatever $\z^1$ in a convex neighborhood of $\z^0$ in $\el^{\infty}_{\zz^-}$ is. If $\z^1\in \el^{\infty}_{\zz^-}$
is arbitrary we can joint it with $\z^0$ by an appropriate polygonal curve and finish the proof in a finite
number of steps.

\smallskip\qed

\newprg[LL.fin]{Lerch-Lipschitz zeta function of finitely many variables}

For $\a=(a_1,...,a_k)\in \cc^k_{\im^+}$ and $\z = (z_1,...,z_k)\in \cc^k_{\re^+}$ we considered in 
the Introduction the following ''cut`` of the Lerch-Lipschitz zeta function $\zeta^{LL}$ by 
functions of finitely many variables
\begin{equation}
\eqqno(ll-f2)
\zeta_k(\a, \z,s) = \sum_{n=1}^{\infty}\frac{e^{2\pi in a_{[n-1]_k+1}}}{(n + z_{[n-1]_k+1})^s} = 
\sum_{p=0}^{\infty}\sum_{j=1}^k
\frac{e^{2\pi i (pk +j) a_j}}{(pk + j + z_j)^s}.
\end{equation}
In more details $\zeta_k(\a , \z, s) $ is the sum of blocks indexed by $p=0,1,...$ . Each 
of these blocks depends on $k$ variables $z_1,...,z_k$ and $k$ parameters $a_1,...,a_k$ as follows
\begin{equation}
\eqqno(blocks)
\zeta_k(a, z,s) = \underbrace{\frac{e^{2\pi i a_1}}{(1+z_1)^s} + ... + 
\frac{e^{2\pi i k a_k}}{(k+z_k)^s}}_{k \text{ terms }} +
... + \underbrace{\frac{e^{2\pi i (pk+1)a_1}}{(pk+1+z_1)^s} + ... + 
\frac{e^{2\pi i (pk+k)a_k}}{(pk+k+z_k)^s}}_{k \text{ terms }}
+ ...
\end{equation}
When $k\to \infty$ functions $\zeta_k(\{a_j\}_{j=1}^k, \{z_j\}_{j=1}^k,s)$ in some sense approximate 
the function $\zeta (\a,\z,s) = \zeta (\{a_j\}_{j=1}^{\infty},\{z_j\}_{j=1}^{\infty},s)$. From \eqqref(blocks)
it is well seen that $\zeta_k$ is well defined for $\z = (z_1,...,z_k)$ with $z_j\in \cc\setminus \{z\in \zz^-: [z]_k=-j\} = \cc\setminus \zz^-_j$ and $s\in H_1$. It is clear as well that $\zeta_k$ is continuous in 
parameter $\a$ up to $\im a_j=0$.

\smallskip Function $\zeta_k$ can be obtained from the Euler-Riemann form of $\zeta$ by the substitution

\begin{equation}
\eqqno(chg-var1)
\begin{cases}
b_n = e^{2\pi i na_{[n-1]_k+1}},\cr 
z_n = \frac{1}{n+z_{[n-1]_k+1}}, 
\end{cases}
\end{equation}
where $\im a_j >0$ and $\re z_j >0$ for $j=1,...,k$, and $n\in \nn $. Therefore it follows also from Corollary \ref{l2-cont1}, which will be proved later, that $\zeta_k$ admits a multivalued analytic continuation to 
$\cc^k_{\im^+}\times \cc^k_{\zz^-}\times H_1$ continuous in parameter $\a$ up to $\bar\cc^k_{\im^+}$. We shall
not stop on this here. 

\newprg[LL.real]{Analytic continuation with respect to variables $\s$ and $\z$}

For a fixed $\a\in \bar\cc^k_{\im^+}$ and fixed $\z \in \cc^k_{\zz^-}$ let us fulfill an analytic continuation 
of $\zeta_k(\a,\z,s)$ in $s$-variable from $H_1$ to $\cc$. First we restrict $\z$-variable to $\rr^k_+$, \ie 
$\z = \x \in \rr^k_+\deff (\rr_+)^k$, where $\rr_+ =\{t\in \rr :t\ge 0\}$.

\begin{figure}[h]
\centering
\includegraphics[width=1.8in]{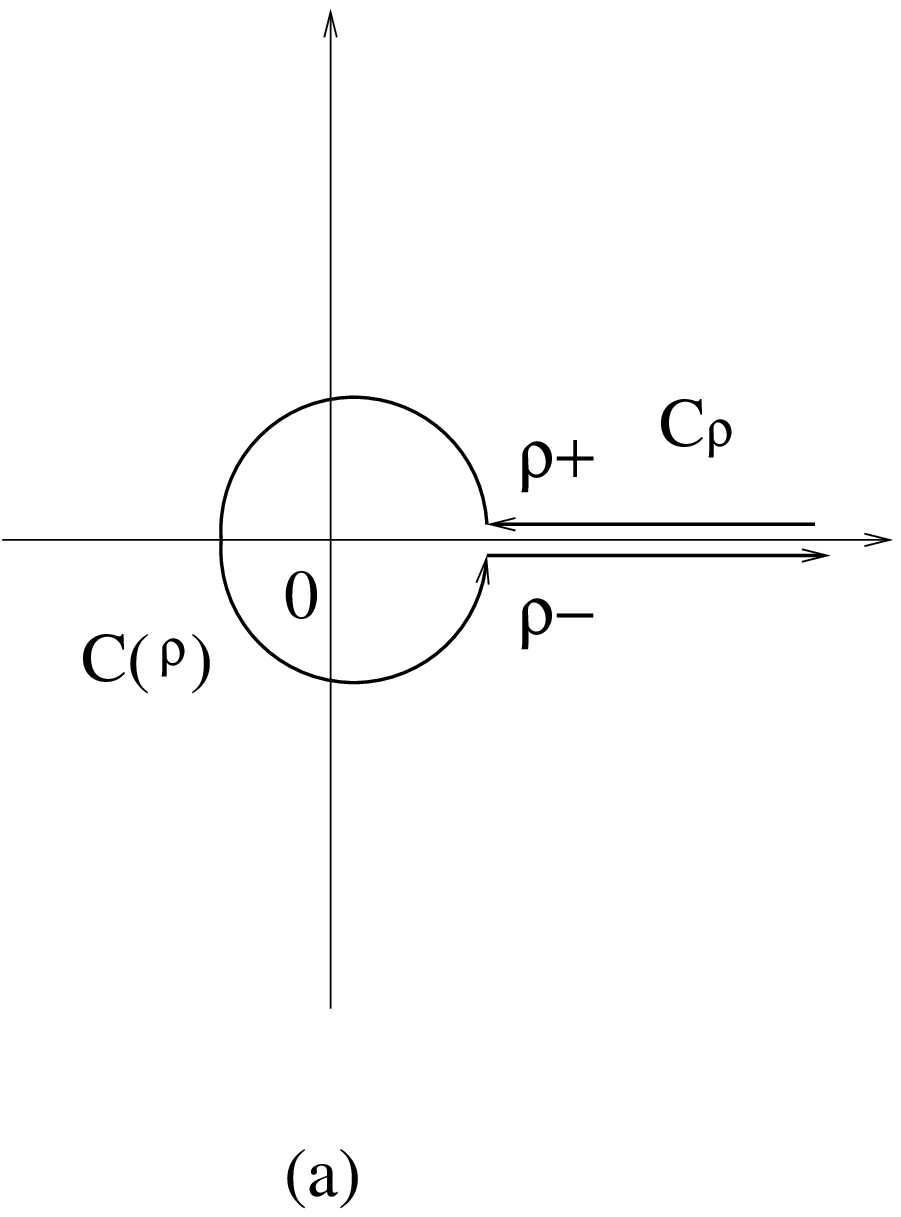}
\quad\quad\quad\qquad
\includegraphics[width=1.9in]{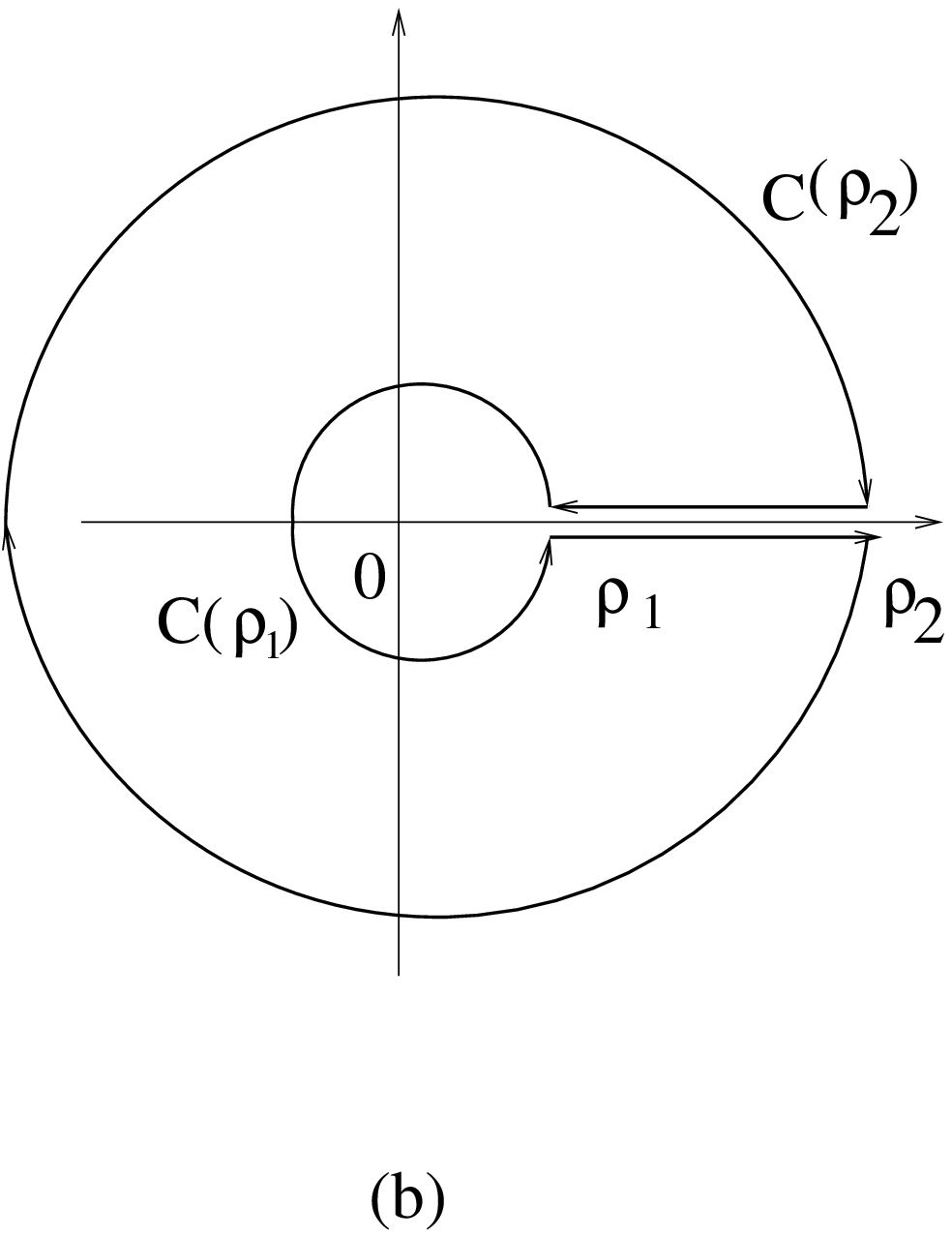}
\caption{{\bf (a)} $C_{\rho} = (+\infty , \rho^+)\cup C(\rho)\cup (\rho^-, +\infty)$, where by $\rho^{\pm}$
we denote the same point $\rho\in \rr^+$ but with different values of the argument: $\arg (\rho^+) = 0$ 
and $\arg (\rho^-) = 2\pi$. $C(\rho)$ is the circle of radius $\rho$ centered at the origin. {\bf (b)} 
Integral along $C_{\rho_1}$ can be rewritten as an integral along $C_{\rho_2}$
plus an integral along the closed contour which is drawn on this picture. This closed contour doesn't bound 
points as in \eqqref(zet-k) and therefore these integrals are equal.}
\label{contour1-fig}
\end{figure}

\begin{lem}
\label{cont-x-s}
For every couple $(\a,\x)\in \bar\cc^k_{\im^+}\times\rr_+^k$ function $\zeta_k(\a,\x,s)$ can be analytically 
continued in $s$ to a meromorphic function on $\cc$. Moreover:

\sli if $a_j\not\in \zz\left[1/k\right]$ for all $j\in \{1,...,k\}$ then $\zeta_k(\a,\x,s)$ is
holomorphic on $\cc$;
        
\slii if $a_j\in \zz\left[1/k\right]$ for some $j\in \{1,...,k\}$ then $\zeta_k(\a,\x,s)$ has a 
simple pole at $s=1$  

\quad with the residue equal to 
\begin{equation}
\eqqno(res-zeta1)
\res (\zeta_k,1) = 1/k \sum_{a_j\in \zz\left[1/k\right]}e^{2\pi ija_j}.
\end{equation}
\end{lem}
\proof The proof follows the classical line of arguments, which are due to Riemann, see \cite{Rm}. Making the substitution 
$t\to (pk + j + x_j)t$ to the Euler's gamma function 
\[
\Gamma (s) = \int\limits_0^{\infty}t^{s-1}e^{-t}dt = (pk + j + x_j)^s \int\limits_0^{\infty}t^{s-1} e^{-(pk + j + x_j)t}dt
\]
we get from \eqqref(ll-f2)
\[
\zeta_k(\a,\x,s) = \frac{1}{\Gamma (s)}\sum_{p=0}^{\infty}\sum_{j=1}^k\int\limits_0^{\infty}t^{s-1}e^{2\pi i (pk +j) a_j}
 e^{-(pk + j + x_j)t}dt = 
\]
\begin{equation}
\eqqno(ll-f3)
 = \frac{1}{\Gamma (s)}\int\limits_0^{\infty}t^{s-1}\sum_{p=0}^{\infty}\sum_{j=1}^ke^{2\pi i (pk +j) a_j}
 e^{-(pk + j + x_j)t}dt
\end{equation}
due to the absolute convergence of this series for $s\in H_1$ and $\im a_j \ge 0$. We can rewrite \eqqref(ll-f3) as 
follows
\[
\zeta_k(\a,\x,s) = \frac{1}{\Gamma (s)}\sum_{j=1}^ke^{2\pi i ja_j}\int\limits_0^{\infty}t^{s-1} e^{-(j + x_j)t}
\sum_{p=0}^{\infty}e^{2\pi i pka_j} e^{-pkt}dt = 
\]
\begin{equation}
\eqqno(ll-f4)
= \frac{1}{\Gamma (s)}\sum_{j=1}^ke^{2\pi i ja_j} \cdot\int\limits_0^{\infty}t^{s-1} e^{-(j + x_j)t}\sum_{p=0}^{\infty}
\left(e^{2\pi i ka_j-kt}\right)^pdt =
\end{equation}
\[
= \frac{1}{\Gamma (s)}\sum_{j=1}^ke^{2\pi i ja_j} \int\limits_0^{\infty}t^{s-1}
\frac{e^{-(j + x_j)t}}{1 - e^{2\pi i ka_j-kt}}dt.
\]
The integral on the right hand side of \eqqref(ll-f4) converges absolutely for $\re s>1$ and 
$(\a , \x)\in \bar\cc^k_{\im^+}\times \rr^k_+$. For every $j=1,...,k$ and small $\rho >0$  consider the integral 
\begin{equation}
\eqqno(i-rho-0)
I^j_{\rho}(a_j,x_j,s) \deff \int\limits_{C_{\rho}}\lambda^{s-1} \frac{e^{-(j + x_j)\lambda}}{1 - e^{2\pi i ka_j - k\lambda }}d\lambda,
\end{equation}
where $\lambda = t+i\tau$ and $\C_{\rho}$ is the contour defined below, see Figure \ref{contour1-fig} (a).
Let us explain that $I^j_{\rho}(s)$ doesn't depend on $\rho >0$. Indeed, the only problem 
could arise with denominator if $2\pi i ka_j - k\lambda \in 2\pi i\zz$, say equals to $2\pi i l$ for some 
$l\in \zz$. This means that 
\begin{equation}
\eqqno(zet-k)
\lambda = 2\pi i \left(a_j - \frac{l}{k} \right).
\end{equation}
In the case under discussion we take $\rho >0$ small enough in order that the contour $C_{\rho}$ avoids points as in 
\eqqref(zet-k) for all $j=1,...,k$ and, moreover, the disk bounded by $C(\rho)$ doesn't contains these points except,
possibly, the origin. 

This implies holomorphicity of $I_{\rho}^j(s)$ on the whole of $\cc$. At the same time for 
$s\in H_1$ one has 
\begin{equation}
\eqqno(i-rho-1)
I^j_{\rho}(a_j,x_j,s) = \lim_{\substack{\rho \to 0}}\int\limits_{C_{\rho}}\lambda^{s-1}\frac{e^{-(j + x_j)\lambda}}{1 -
e^{2\pi i ka_j - k\lambda}}d\lambda =(e^{2\pi i s} - 1) \int\limits_0^{\infty}t^{s-1} 
\frac{e^{-(j + x_j)t}}{1 - e^{2\pi i ka_j - kt}}dt.
\end{equation}
Therefore from \eqqref(ll-f4) we see that for $s\in H_1$ one has
\begin{equation}
\eqqno(ll-f5)
\left(e^{2\pi i s } - 1\right)\Gamma (s)\zeta_k(\a,\x,s) = \sum_{j=1}^ke^{2\pi i ja_j}I^j_{\rho}(a_j,x_j,s),
\end{equation}
with the right hand side being holomorphic on $\cc$. This gives the analytic continuation of $\zeta_k(\a,\x,s)$ to a 
meromorphic function on $\cc$ with eventual poles at points of $\zz_{>0}\deff \{n\in \zz : n>0\}$ at most.

\smallskip\noindent{\slsf Case 1.} {\it $a_j\not\in \zz \left[1/k\right]$ for all $j\in \{1,...,k\}$.} In this 
case $\zeta_k(\a,\x,s)$ is an entire function. Indeed, observe that for $s\in \zz$ one has 
\begin{equation}
\eqqno(i-rho5)
I^j_{\rho}(a_j,x_j,s) = \int\limits_{C_{\rho}}\lambda^{s-1} \frac{e^{-(j + x_j)\lambda}}{1 - e^{2\pi i ka_j - k\lambda}}
d\lambda = \int\limits_{C(\rho )}\lambda^{s-1} \frac{e^{-(j + x_j)\lambda}}{1 - e^{2\pi i ka_j - k\lambda}}d\lambda = 
\end{equation}
\[
=
\begin{cases}
 0 \text{ if } s\in \zz\setminus \{0\},\cr
 \frac{2\pi i }{1- e^{2\pi ika_j}}\text{ if } s=0. 
\end{cases}
\]
This implies that $\zeta_k(\a,\x,s)$ is an entire function which vanishes when $s$ is a negative  integer, 
\ie $s\in \zz_{<0}\deff \{-1, -2,...\}$. It satisfies the following relation on the whole of $\cc$
\begin{equation}
\eqqno(ll-f6)
\zeta_k(\a,\x,s) = \frac{1}{(e^{2\pi i s} -1)\Gamma (s) } \sum_{j=1}^ke^{2\pi i ja_j}I_{\rho}^j(a_j,x_j,s).
\end{equation}

\smallskip\noindent{\slsf Case 2.} { $a_j\in \zz\left[1/k\right]$ for some $j\in \{1,...,k\}$.} In that
case for such $j$ and $s\in \zz$ we have
\begin{equation}
\eqqno(i-rho4)
I^j_{\rho}(a_j,x_j,s) \deff \int\limits_{C_{\rho}}\lambda^{s-1} \frac{e^{-(j + x_j)\lambda}}{1 - e^{ - k\lambda}}d\lambda =
\int\limits_{C(\rho)}\lambda^{s-1} \frac{e^{-(j + x_j)\lambda}}{1 - e^{ - k\lambda}}d\lambda  = 
\begin{cases}
0 \text{ if } s\in \zz \setminus \{1\},\cr
\frac{2\pi i }{k}\text{ if } s=1. 
\end{cases}
\end{equation}
The same relation \eqqref(ll-f5) or, better \eqqref(ll-f6), provides us the meromorphic continuation of $\zeta_k(\a,\x,s)$
to $\cc$ having a simple pole at $s=1$ with residue 
\[
\res (\zeta_k,1) = \lim_{s\to 1}\left[\frac{s-1}{e^{2\pi is} -1}\frac{1}{\Gamma (s)}\sum_{j=1}^k
e^{2\pi ija_j}I_{\rho}^j(s)\right] = 1/k\sum_{a_j\in \zz\left[1/k\right]}e^{2\pi ija_j}.
\]

\smallskip\qed

\smallskip Recall that a subset $R$ of an open set $D$ in a complex manifold called {\slsf plurithick} if 
for every sequence $\{u_n\}$ of uniformly bounded from above plurisubharmonic functions on $D$ condition 
\begin{equation}
\eqqno(p-thick-1)
\limsup_{\substack{n\to \infty}}u_k|_R \equiv - \infty \quad\text{ implies } \quad 
\limsup_{\substack{k\to \infty}}u_k \equiv - \infty .
\end{equation}

\begin{lem}
\label{cont-z-s}
For every  $\a\in \bar\cc^k_{\im^+}$ function $\zeta_k(\a,\cdot,\cdot)$ can be 
analytically continued along any path in $\cc^k_{\zz^-}\times \cc$. Moreover:

\sli if $a_j\not\in \zz\left[1/k\right]$ for all $j\in \{1,...,k\}$ then $\zeta_k(\a,\z,\cdot)$ is holomorphic 
on $\cc$;
        
\slii if $a_j\in \zz\left[1/k\right]$ for some $j\in \{1,...,k\}$ then $\zeta_k(\a,\z,\cdot)$ has a simple pole 
at $1$ with the

\quad  residue equal to 
\begin{equation}
\eqqno(res-zeta2)
\res (\zeta_k,1) = 1/k\cdot \sum_{a_j\in \zz\left[1/k\right]}e^{2\pi ija_j}.
\end{equation}

\sliii In addition for every $j=1,...,k$ and every $p\in \zz^-$ the monodromy of $\zeta_k(\a , \cdot , s)$ around 

\quad the hyperplane $z_j = - pk - j$ is equal to $\frac{e^{2\pi i (pk+j)a_j}}{(pk+j+z_j)^s}\left(e^{-2\pi i s}-1\right)$.
\end{lem}
\proof (i) Take any path $\gamma : [0,1]\to \cc^k_{\zz^-}$ starting at $\adyn =(1,...,1)\in \rr^k_+$. Let $B_t$ 
be the ball $B_t(\gamma (t),\eps)$ in $\cc^k$ with an appropriate $\eps >0$. We need to extend $\zeta_k(\a , \cdot ,\cdot)$
to $B_t\times \cc$ for all $t$. Since $\rr^k_+\cap B_0$ is a plurithick subset of $B_0$, as well as $B_{t_1}\cap B_t$
in $B_t$ for every $t_1<t$ close to $t$ the needed statement follows from Corollary \ref{siciak-c}. 

\smallskip\noindent (ii) Repeat the same agrument for the function $(s-1)\zeta_k(\a , \cdot , \cdot)$ to get its
holomorphic extension. This provides the meromorphic extension of $\zeta_k$ itself with the simple pole at $\adyn$. As
for the residu remark that 
\[
\res (\zeta_k(\a , z, \cdot) , \adyn) = \frac{1}{2\pi i}\int\limits_{|s-1| = \eps}\zeta_k(\a , z,s)ds
\]
holomorphically depends on $z$ and therefore is constant, because it is constant on $\rr^k_+$.

\smallskip\noindent (iii) The monodromy of $\zeta_k(\a , \cdot , s)$ around the hyperplane $z_j = -pk - j$ 
in the positive direction is, in fact, the monodromy  along the closed path in $z_j$-plane starting at $z_j\not=-pk-j$ 
and going around $-pk-j$ in the counterclockwise direction of the summand
\begin{equation}
\eqqno(monod-0)
\frac{e^{2\pi i (pk+j)a_j}}{(pk+j+z_j)^s},
\end{equation}
which is obviously equal to $\frac{e^{2\pi i (pk+j)a_j}}{(pk+j+z_j)^s}\left(e^{-2\pi i s}-1\right)$. 

\smallskip\qed

\begin{rema} \rm
\label{monod-z}
Notice that item (\sliii of this lemma implies that $\zeta (\a , \cdot , s)$ descends to the 
abelian cover of $\cc^k_{\zz^-}\times \cc$, \ie to $\tilde\cc^{k, ab}_{\zz^-}\times \cc$ and 
separates points over $\cc^k_{\zz^-}\times \cc$ as stated in Theorem \ref{ll-thm}. Indeed, 
the fundamental group of $\cc_{\zz^-}^k$ is generated by simple loops starting from $\adyn = 
(1,...,1)$ and going around points $-pk-j$ on the $z_j$-plane for all $p\in \nn$ and $1\le j\le k$. 
Let $\gamma_1 $ be such loop around $-p_1k-j$ on the $z_j$-plane and $\gamma_2$ around $-p_2k-r$ 
on the $z_r$-plane. Then the monodromy of $\zeta_k(\a , \cdot , s)$ along $\gamma_1^{-1}\cdot 
\gamma_2^{-1}\cdot\gamma_1\cdot\gamma_2$ is obviously zero if $r\not=j$. If $r=j$ then it is 
zero too. Therefore the monodromy along any loop generating the commutator of $\pi_1(\cc_{\zz^-}^k)$
is trivial. 
\end{rema}

\newprg[LL.a]{Analytic continuation with respect to the parameter $\a$}

Now let us turn our attention to the parameter $\a$. Up to now it was aloud to vary in $\bar\cc^k_{\im^+}$.
We shall prove  the following

\begin{lem}
\label{cont-a}
Function $\zeta_k$ admits an analytic continuation to a meromorphic function on 
\begin{equation}
\tilde\calz_k \deff \widetilde{\cc}^{k,ab}_{\zz [1/k]}\times \widetilde{\cc}^{k,ab}_{\zz^-}\times \cc ,
\end{equation}
where $\widetilde{\cc}^{k,ab}_{\zz[1/k]}$ is the abelian cover of $\cc^k_{\zz[1/k]} \deff
\left(\cc\setminus \zz \left[1/k\right]\right)^k$. The monodromy of $\zeta_k(\cdot , \z ,s)$ 
around the hyperplane $\{a_j=l/k\}$  is equal to
\begin{equation}
\eqqno(monod-aj)
\frac{[2\pi (a_j-\frac{l}{k})]^{s-1}e^{i\frac{\pi}{2}(s-1)}e^{-2\pi i(a_j-\frac{l}{k})(j+z_j)}}{k}.
\end{equation}

\end{lem}
\proof This issue will be proved by studying integrals from the right hand part of \eqqref(i-rho-1). 
We follow \cite{LL,Le}. For a fixed $s\in H_1$ and $x\in \rr^k_+$ write as in \eqqref(i-rho-1) 
$I_p^j(a_j, x_j,s) = (e^{2\pi i s}-1)J_p^j(a_j,x_j,s)$, \ie 

\begin{equation}
\eqqno(j-rho1)
J_{\rho}^j(a_j,x_j,s)\deff \int\limits_0^{\infty}t^{s-1}\frac{e^{-(j + x_j)t}}{1 - e^{2\pi i ka_j - kt}}dt.
\end{equation}
One readily sees that this integral is well defined and holomorphic with respect to $a_j$ unless 
\begin{equation}
\eqqno(aj-1)
a_j-\frac{l}{k}\in -\frac{i}{2\pi}\rr^+ \quad \text{ for some } l\in\zz .
\end{equation}
\begin{figure}[h]
\centering
\includegraphics[width=2.6in]{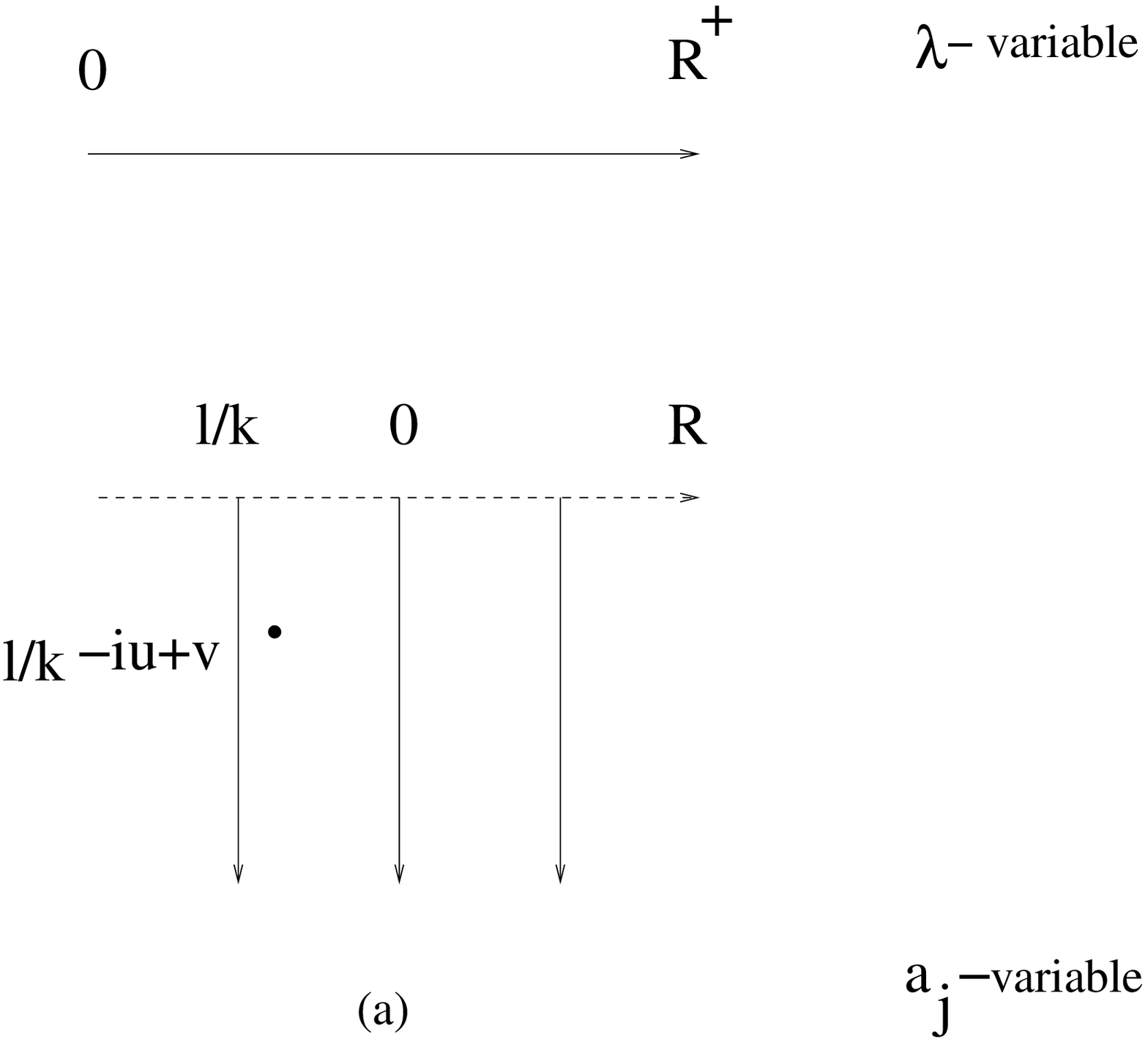}
\quad\quad\quad\qquad
\includegraphics[width=1.5in]{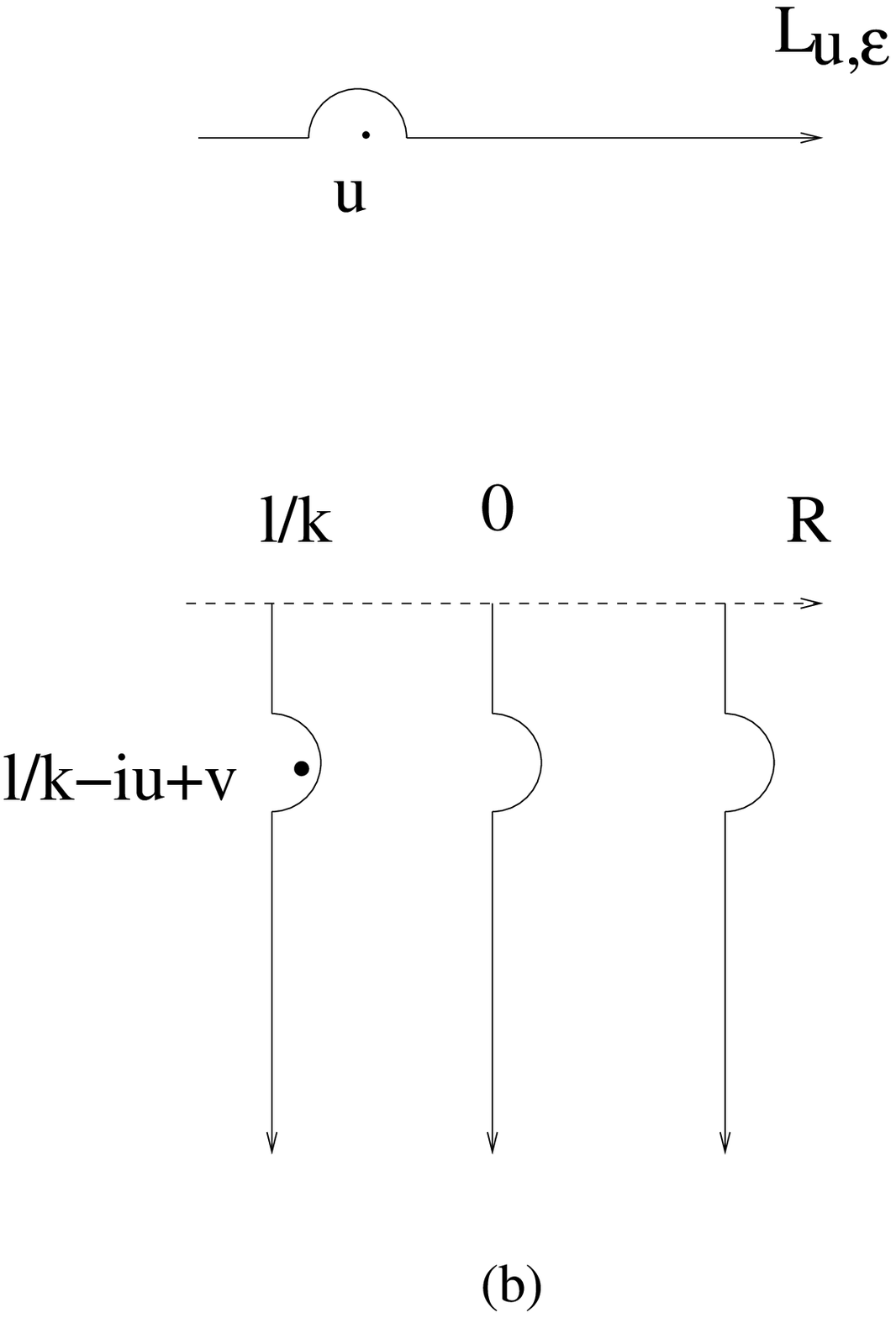}
\caption{{\bf (a)} Integration along $\rr^+$ gives us zeta function defined for $a_j$ on $\cc$ minus 
vertical halph-lines starting at points of $\zz\left[1/k\right]$. {\bf (b)} Integration along
$L_{u,\eps}$ produces extension to bumps as on the right. To find monodromy of $\zeta_k$ one should take 
the difference between the integral along $\rr^+$ and $L_{u,\eps}$ for the values of $a_j=\frac{l}{k}-iu+v$, 
which is marked by a bold point on the both Pictures (a) and (b).}
\label{contour2-fig}
\end{figure}
In other words $J_{\rho}^j(a_j,x_j,s)$ is holomorphic in $a_j$ on $\cc\setminus \left\{\zz\left[1/k\right]
-i\rr^+\right\}$. Fix some $u\in \rr^+$ and deform $\rr^+$ to a contour $L_{u,\eps} = [0,u-\eps]\cup S_{\eps}^+(u)
\cup [u+\eps, +\infty)$ as on the Figure \ref{contour2-fig}. Here $S_{\eps}^+(u)$ is the upper part of the 
circle of radius $0<\eps <u/2$ centered at $u$. We have that 
\begin{equation}
\eqqno(j-rho2)
\int\limits_0^{\infty}t^{s-1}\frac{e^{-(j + x_j)t}}{1 - e^{2\pi i ka_j - kt}}dt = 
\int\limits_{L_{u,\eps}}\lambda^{s-1}\frac{e^{-(j + x_j)\lambda}}{1 - e^{2\pi i ka_j - k\lambda}}d\lambda .
\end{equation}
Indeed, for $\lambda = u + \eps e^{i\theta}, 0\le \theta \le \pi$, the condition  \eqqref(zet-k) reads as 
\begin{equation}
\eqqno(zet-k1)
a_j-\frac{l}{k} = -\frac{i(u+\eps\cos \theta)}{2\pi} + \frac{\eps\sin \theta}{2\pi},
\end{equation}
and this is unconsistent for $\im a_j \ge 0$ and the choice of $\eps$ made. Deforming the contour of 
integration from $\rr^+$ to 
$L_{u,\eps}$ we extend $J_{\rho}^j(s)$ analytically in the variable $a_j$ to an $\eps$-neighborhood
of each point $\frac{l}{k} - iu$ from the left, see Figure \ref{contour2-fig}(b), \ie to a point of 
the form $l/k -iu+v$ with $v>0$. The monodromy calculated at such point $a_j=\frac{l}{k}-iu+v$ with 
$v>0$ is equal to the difference between the left-hand and the right-hand sides of \eqqref(j-rho2).
Taking into account that $u+iv = i(a_j-l/k)$ this monodromy is equal to 
\ie to
\[
\int\limits_{S^+_{\eps}(u)+[u-\eps , u+\eps]} \lambda^{s-1}\frac{e^{-(j+x_j)\lambda}}
{1-e^{2\pi i k (\frac{l}{k}-iu+v)-k\lambda}}d\lambda =
\int\limits_{S^+_{\eps}(u)+[u-\eps , u+\eps]} \lambda^{s-1}\frac{e^{-(j+x_j)\lambda}}
{1-e^{2\pi k (u+iv-\lambda/2\pi)}}d\lambda =
\]
\[
=\lim_{\lambda \to 2\pi i(a_j-\frac{l}{k})}(\lambda - 2\pi i(a_j-\frac{l}{k})) 
\lambda^{s-1}\frac{e^{-(j+x_j)\lambda}}{1-e^{k(2\pi i(a_j-\frac{l}{k})-\lambda)}} =
\]
\begin{equation}
\eqqno(monod-2)
=\frac{[2\pi i(a_j-\frac{l}{k})]^{s-1}e^{-2\pi i(a_j-\frac{l}{k})(j+x_j)}}{k}
= \frac{[2\pi (a_j-\frac{l}{k})]^{s-1}e^{i\frac{\pi}{2}(s-1)}e^{-2\pi i(a_j-\frac{l}{k})(j+x_j)}}{k}.
\end{equation}
Here $\arg i(a_j-\frac{l}{k}) = \arg (u+iv)$ should be taken in $[0,\pi/2)$, and therefore 
$\arg (a_j-\frac{l}{k})
\in [-\pi/2, 0)$. Therefore for $(x, s)\in \rr^k_+\times H_1$ function 
\[
\left(e^{2\pi i s } - 1\right)\Gamma (s)\zeta^{LL}_k(a ,x,s) = \sum_{j=1}^ke^{2\pi i ja_j}I^j_{\rho}(s)
\]
extends in variable $a$ analytically from $\cc^k_{\im^+}$ to the abelian cover 
$\widetilde{\cc}^{k,ab}_{\zz[1/k]}$ of $\cc^k_{\zz[1/k]} = $ $=\left(\cc\setminus \zz \left[1/k\right]\right)^k$. Since $\rr^k_+\times H_1$ is plurithick in $\widetilde{\cc}^{k,ab}_{\zz^-}\times \cc$ we obtain that the same is true for all $(z ,s)\in \widetilde{\cc}^{k,ab}_{\zz^-}\times \cc$. Lemma and Theorem \ref{ll-thm} are
proved.

\smallskip\qed

\newappx[A1]{Pluricomplex analysis on complex sequence spaces}

We recall here some rudiments of pluricomplex analysis on sequence spaces which we used along the text. 
We need to do this since the Hartogs-Kazaryan type results we used for the proof of the principal results 
of this paper seem not to be redily avaiable in the literature in the infinite dimensional case. Our
exposition will cover finite and infinite diemnsinal cases simultaneously.

\newprg[APP1.hol]{Holomorphic functions}

We shortly recall a few needed notions concerning Banach holomorphicity in order to fix the notations adapted to 
our sequence spaces. For details we refer to  \cite{Mu}. For a natural $m\ge 1$  an $m$-linear form on $\el_{\r}$ 
is a  mapping 
\[
A:\underbrace{\el_{\r}\times ... \times \el_{\r}}_{m-times}\to \cc ,
\]
which is linear on each variable. Mapping $A$ is called $\r$-bounded if 
\begin{equation}
\eqqno(norm-am)
\norm{A}_{\r} \deff \sup\{|A(z_1,...,z_m)|: z_j\in \bar B_{\r}^{\infty}(0,1)\} <\infty .
\end{equation}
We denote by $\call_m^{\b}(\el_{\r})$ the linear space of $\r$-bounded $m$-linear forms on $\el_{\r}$. 

%\underbrace{\el_{\r}\times ... \times \el_{\r}}_{m-\text{times}}%

\begin{defi}
An $\r$-bounded homogeneous polynomial of degree $m$ on $\el_{\r}$ is such mapping $P:\el_{\r}\to \cc$ that
\begin{equation}
\eqqno(m-pol)
P(z) = A(\underbrace{z,...,z}_{m-\text{times}})
\end{equation}
for some $A\in \call_m^{\b}(\el_{\r})$.
\end{defi}

An $m$-linear form $A$ defining $P$ can be supposed to be symmetric because its {\slsf symmetrization}
\[
A^{\sym}(z_1,...,z_m) = \frac{1}{m!}\sum\limits_{\sigma\in S_m}
 A(z_{\sigma (1)},...,z_{\sigma (m)}). 
\]
defines the same polynomial $P$ and, moreover, $A^{\sym}$ is $\r$-bounded if such is $A$. From the polarization 
formula for a symmetric $A$, \ie from  
 \begin{equation}
\eqqno(polariz) A(z_1,...,z_m) =
\frac{1}{m!2^m}\sum\limits_{\eps_j=\pm
1}\eps_1...\eps_mA(z_0+\eps_1z_1+...+\eps_mz_m,...,z_0+\eps_1z_1+...+\eps_mz_m)
\end{equation}
for every couple $z_0,...,z_m$ of vectors from $\el_{\r}$, see \cite{Mu}, one obtains that such $A$ is uniquely 
defined by $P$. This uniquely for $P$ defined symmetric form let us denote as $\hat P$. Again from the polarization
formula one can see that 
\begin{equation}
\eqqno(hat)
\norm{P}_{\r}\le \norm{\hat P}_{\r} \le \frac{m^m}{m!}\norm{P}_{\r},
\end{equation}
where $\norm{\hat P}_{\r}$ is defined as in \eqqref(norm-am) and the norm of a polynomial $P$ is defined as 
\begin{equation}
\eqqno(norm-pm)
\norm{P}_{\r} \deff \sup\{|P(z)|: z\in \bar B_{\r}^{\infty}(0,1)\}.
\end{equation}

\begin{rema} \rm 
\label{pol-cont}
Notice that an $\r$-bounded $m$-homogeneous polynomial $P$ is {\slsf bounded} and $\r$-continu\-ous on 
$B_{\r}^{\infty}(0, \beta )$ for any $0< \beta < +\infty$. Boundedness is simple and follows from the relation 
\begin{equation}
\eqqno(m-hom)
P(\lambda z) = \lambda^m P(z),
\end{equation}
\ie on $B^{\infty}_{\r}(0,\beta )$ our polynomial $P$ is bounded by $\beta^m\norm{P}_{\r}$. As for the continuity 
take some $z^1,...,z^m,w^1,...,w^m\in B_{\r}^{\infty}(0, \beta)$ such that $z^j-w^j\in B_{\r}^{\infty}(0,\eps)$ for 
$j=1,...,m$. Then 
\[
\left|\hat P(z^1,...,z^m) - \hat P(w^1,...,w^m)\right|\le \left|\hat P(z^1,...,z^m)-\hat P(w^1,z^2,...,z^m)\right|
+ ... +
\]
\[
\left|\hat P(w^1,...,w^{m-1},z^m) - \hat P(w^1,...,w^m)\right| \le \sum_{j=1}^m\left|\hat P(w^1,...,w^j,z^j-w^j, 
z^{j+1},...,z^m)\right|
\le
\]
\[
\le m\eps \beta^{m-1}\norm{\hat P}_{\r},
\]
which proves the assertion.
\end{rema}

\begin{defi}
\label{d-hol-f1}
Let $D\subset \el_{\r}$ be an $\r$-open set. Function $f: D\to \cc$ is called $\r$-holomorphic or, simply holomorphic 
if $\r$ is clear from the context, if the following holds.

\sli $f$ is $\r$-bounded, \ie for every $\beta >0$ such that $\bar B_{\r}^{\infty} (z^0,\beta)\subset D$
function $f$ is 

\quad bounded on $\bar B_{\r}^{\infty} (z^0,\beta_1)$ for any $0<\beta_1<\beta$.

\slii $f$ is G\^ateaux differentiable, \ie  for  every $z\in \el_{\r}$ the function of one complex

\quad   variable $f(z^0+tz)$ is holomorphic in $t$ near the origin.
\end{defi}

\begin{prop}
\label{d-hol-f2}
Function $f$ defined on an $\r$-open set $D$ is $\r$-holomorphic if and only
if for every $\bar B_{\r}^{\infty}(z^0,\beta)\subset D$ there exists  a sequence of $\r$-bounded  
$m$-homogeneous  polynomials $\{P_m\}_{m=1}^{\infty}$ such that 
\begin{equation}
\eqqno(taylor1)
f(z^0 + z) = f(z^0) + \sum_{m=1}^{\infty}P_m(z)
\end{equation}
normally on $\bar B_{\r}^{\infty} (0,\beta_1)$ for every $0<\beta_1<\beta$. The latter means that 
\begin{equation}
\eqqno(taylorm)
\sum_{m=1}^{\infty}\beta_1^m\norm{P_m}_{\r} < \infty .
\end{equation}
\end{prop}
\proof Polynomials $P_m$ depend also on $z^0$ of course and it would be more colloquial to write \eqqref(taylor1) as 
\begin{equation}
\eqqno(taylorz0)
f(z^0 + z) = f(z^0) + \sum_{m=1}^{\infty}P_m(z^0,z).
\end{equation}

\smallskip\noindent $\Leftarrow$ Take $z^0\in D$ and $\beta >0$ such that $\bar B_{\r}^{\infty}(z^0,\beta)\subset D$. 
Due to the assumption there exist $\r$-bounded $m$-homogeneous  polynomials $P_m$ such that 
\eqqref(taylor1) and \eqqref(taylorm) hold on $\bar B^{\infty}(0,\beta_1)$ whatever $0<\beta_1<\beta$ is. This immediately 
implies that $f$ is bounded on all such $\bar B_{\r}^{\infty}(z^0,\beta_1)$. As for the G\^ateaux 
differentiability of $f$ take any $z\in B_{\r}^{\infty}(0,\beta_1)$  and write
\begin{equation}
\eqqno(taylor-t)
f(z^0 + tz) = f(z^0) + \sum_{m=1}^{\infty}P_m(tz) = f(z^0) + \sum_{m=1}^{\infty}P_m(z)t^m.
\end{equation}
Series \eqqref(taylor-t) converge uniformly on $|t| \le 1$ due to \eqqref(taylorm), thus providing the holomorphicity 
of $f(z^0 + tz)$ near the origin.

\smallskip\noindent $\Rightarrow$ Let again $\bar B^{\infty}_{\r}(z^0,\beta)\subset D$. Fix $0<\beta_1<\beta$.
For $z\in B^{\infty}_{\r}(0,\beta_1)$ denote by $\Delta_zf(z^0)$ the derivative at $t=0$ of $f(z^0 + tz)$. 
Take some $z^1,...,z^m\in B^{\infty}_{\r}(0,\beta_1)$. Then for any couple $(t_1,...,t_m)$ in a neighborhood $U$ 
of zero of $\cc^m$ vector $t_1z^1+...+t_mz^m$ belongs to $B^{\infty}_{\r}(0,\beta_1)$. Neighborhood $U$ 
being fixed for a choice of $z^1,...,z^m$ made. By the finite dimensional Hartogs' separate analyticity 
theorem function $f(z^0+t_1z^1+...+t_mz^m)$ is holomorphic in $U$. Therefore we can set 
\begin{equation}
\eqqno(a-m)
A_m(z_1,...,z_m) \deff \Delta_{z^1}\circ ...\circ \Delta_{z^m}f(z^0) = 
\frac{\d^m f(z^0 + t_1z^1+...+t_mz^m)}{\d t_1...\d t_m}\Big|_{t_1=...=t_m=0}.
\end{equation}
Function $A_m$ is obviously symmetric and $m$-linear, this is a finite dimensional argument. Also from one 
complex variable applied to the holomorphic function $f(z^0 + tz)$, where $z\in B^{\infty}_{\r}(0,\beta_1)$ 
and $|t|\le \frac{\beta}{\beta_1}$, we get that
\begin{equation}
f(z^0 + z) = f(z^0) + \sum_{m=1}^{\infty}\frac{1}{m!}A_m(\underbrace{z,...,z}_{m-\text{times}}),
\end{equation}
which is nothing else but the Taylor expansion of $f(z^0 + tz)$ at $t=1$. Set $P_m(z) = \frac{1}{m!}A_m(z,...,z)$.
Since $f$ is bounded, say by $M_{\beta_2}(f)$, on $\bar B^{\infty}_{\r}(z^0,\beta_2)$ for some $\beta_1<\beta_2<\beta$
we see from Cauchy formula and the fact that $f(z^0+tz)$ is defined for $|t|\le \frac{\beta_2}{\beta_1}$ that 
\begin{equation}
\eqqno(a-m-bound)
\left|A_m(z,...,z)\right| = \left|\frac{\d^m f(z^0+tz)}{\d t^m}\Big|_{t=0}\right|\le \frac{m!\beta_1^m}{\beta_2^m}
\times M_{\beta_2}(f).
\end{equation}
Therefore $P_m(z)$ is bounded by $(\beta_1/\beta_2)^mM_{\beta}$ on $\bar B^{\infty}_{\r}(0,\beta_1)$, note that 
$\beta_1/\beta_2 <1$. This provides $\r$-boundedness of $P_m$ and the normal on $\bar B^{\infty}_{\r}(0,\beta_1)$
convergence of the series 
\begin{equation}
\eqqno(taylor2)
\sum_{m=1}^{\infty}P_m(z), \quad\text{ because }\quad \sum_{m=1}^{\infty}|P_m(z)|\le M_{\beta} \sum_{m=1}^{\infty}
\left(\frac{\beta_1}{\beta_2}\right)^m = M_{\beta_2}(f) \frac{\beta_1/\beta_2}{1-\beta_1/\beta_2}.
\end{equation}

\smallskip\qed

In the following propositions we list the elementary properties of $\r$-holomorphic functions.

\begin{prop}
\label{unif}
One has the following:

\smallskip \sli A uniform on bounded open sets limit of $\r$-holomorphic functions is $\r$-holomorphic.

\smallskip\slii Let $f$ be $\r$-holomorphic on an $\r$-open  $D\subset \el_{\r}$ and let $\r_1$ be a polyradius
lexicografically 

smaller than $\r$. Then $D_1\deff D\cap \el_{\r_1}$ is $\r_1$-open and $f|_{D_1}$ is $\r_1$-holomorphic.
\end{prop}
\proof \sli The question being local it is sufficient to consider two balls $\bar B^{\infty}_{\r}(0,\beta_1)
\subset B^{\infty}(0,\beta_2)$, \ie $0<\beta_1<\beta_2$ and a sequence $\{f_n\}$ of holomorphic on 
$B^{\infty}(0,\beta_2)$ function uniformly on $B^{\infty}(0,\beta_2)$ converging to $f$. $\r$-boundedness
of the limit $f$ is obvious. The  G\^ateaux holomorphicity of $f$ follows from the Weierstrass theorem.

\smallskip\slii When saying that $\r_1=\{r_{1,n}\}$ is smaller than $\r=\{r_n\}$ lexicografically we mean that 
$r_{1,n}\le r_n$ for all $n$. Remark that $B^{\infty}_{\r_1}(0,1) \subset B^{\infty}_{\r}(0,1)$ in this
case. And this implies the $\r_1$-openess of $D_1$. Continuity of $f|_{D_1}$ follows from the fact that 
$\norm{\cdot}_{\r}\le \norm{\cdot}_{\r_1}$. Its G\^ateaux holomorphicity is obviously follows from that of $f$.

\smallskip\qed

\begin{prop}
\label{uniq}
An $\r$-holomorphic function on an $\r$-open connected set $D$ possesses the following properties.

\smallskip\sli If $f$ vanishes on some $\r$-open subset of $D$ then $f$ vanishes identially. If $|f|$ achievs

\quad its  maximum at some interior point of $D$ then $f\equiv \const$.

\slii $\r$-holomorphic functions  possess the following uniqueness property: let $z^0 = x^0$

\quad $+iy^0 \in D$ and $\beta >0$ be such that $\bar B^{\infty}_{\r}(z^0,\beta )\subset D$. If the 
       restriction of $f$ to $iy^0+$ 

\quad $\bar B^{\infty}_{{\r}, \rr}(x^0, \beta\r)$  vanishes then $f$ vanishes identically.
\end{prop}
\proof Item (\sli can be proved by considering restrictions $f(z^0+tz)$ for $z\in B^{\infty}_{\r}(0,\beta_1)$.

\smallskip\noindent (\slii First we shall prove that $f$
vanishes on $\bar B^{\infty}_{\r}(z^0,\beta )$. Without loss of generality we can assume that $z^0=0$. Write 
the Taylor expansion \eqqref(taylor1) for the real values of variables
\begin{equation}
\eqqno(taylor3)
f(x) = \sum_{m=1}^{\infty}P_m(x),
\end{equation}
which converges on $\bar B^{\infty}_{\r}(x^0, \beta_2)$ for any $0<\beta_1<\beta_2<\beta$. Since for real $t$ one has
that $0\equiv f(tx) = \sum_{m=1}^{\infty}P_m(x)t^m$ we deduce that $P_m(x)=0$ for all $x\in B^{\infty}_{{\r},\rr}
(0,\beta_2)$. Let us prove that $P_m(z)=0$ also for complex $z\in B^{\infty}_{\r}(0,\beta_1)$. Fix such
$z= x+iy$ and remark that $x\in  B_{\r ,\rr}^{\infty}(0,\beta_1\r)$. Therefore the function of one complex 
variable $p(\lambda) \deff P_m(\lambda x)$, which is well defined for $|\lambda|\le \beta_1/\beta_2$, 
vanishes for all $|\lambda|\le \beta_1/\beta_2$. In particular it vanishes at $\lambda_0$ which is such
that $\lambda_0 x = z$. We proved that $P_m$ vanishes on $B_{\r}^{\infty}(0,\beta_1)$. This means that 
$P_m(z)\equiv 0$ in $ B_{\r}^{\infty}(0,\beta\r)$ for all $\beta >0$. Therefore $f$ vanishes on every 
$B_{\r}^{\infty}(0,\beta_1 )$ such that  $0<\beta_1 <\beta$ and $\bar B^{\infty}_{\r}(0,\beta)\subset D$.
The uniqueness property (\sli implies the needed statement.

\smallskip\qed 

\newprg[APP1.sub]{Plurisubharmonic functions}

Let $D$ be a $\r$-open subset in $\el_{\r}$.

\begin{defi}
\label{d-psh}
Function $u:D\to \rr\cup\{-\infty\}$ is called $\r$-plurisubharmonic, $\r$-psh for short, if it is 
upper-semicontinuous in $\r$-topology and for every $z^0\in D$ and every $z\in\el_{\r}$ the restriction 
$u|_{\Omega}$ of $u$ to  $\Omega \deff D\cap \{z^0 + tz:t\in \cc\}$  is subharmonic.
\end{defi}
If $\r$ is clear from the context we say simply that $u$ is psh. The set of $\r$-psh functions on $D$ we shall
denote as $\calp_{\r}(D)$. 

\begin{defi}
An $\r$-open set we call $\calp_{\r}$-hyperconvex if there exists a negative $u\in\calp_{\r}(D)$ such that 
$u(z) \to 0$ when $z \to \d D$ in $\r$-topology.
\end{defi}
Negative here means that $u$ is strictly $< 0$ on the $\r$-interior of $D$.

\begin{exmp}
\label{reg-set1}
The ball $B^{\infty}_{\r}(z^0,\beta)$ is $\calp_{\r}$-hyperconvex for any $\beta >0$. \rm It is sufficient to 
prove this for $\beta =1$. To see this let us prove that the function
\begin{equation}
\eqqno(d-psh-f1)
u(z) \deff \ln\left(\sup_{\substack{n\in \nn}}\left\{\frac{|z_n-z^0_n|}{r_n}\right\}\right) = 
\ln\norm{z-z_0}_{\r},  
\end{equation}
is $\r$-psh on $B^{\infty}_{\r}(z^0,1)$. $u$ is an increasing limit of plurisubharmonic functions of 
finitely many variables
\[
u_N(z) = \ln \left(\max\left\{\frac{|z_n-z^0_n|}{r_n}\right\}_{n=1}^N\right) = 
\max\left\{\ln\frac{|z_n-z^0_n|}{r_n}\right\}_{n=1}^N.
\]
Therefore it is sufficient to prove that it is upper-semicontinuous. For this it is sufficient
to prove the upper-semicontinuity of $u_0(z) \deff e^{u(z)}=\norm{z-z_0}_{\r}$. But the letter is 
continuous in fact.
\end{exmp}

\begin{rema} \rm 
\label{phs-f1}
Remark that function $v_0=u_0-1$ is $\r$-psh on $B^{\infty}_{\r}(z^0,1)$, non-positive on the interior,
$\ge -1$ everywhere, and $\equiv 0$ on the boundary.
\end{rema}

For an $\r$-open set $D$ denote by $\calp_{\r}^-(D)$ the set of non-posive $u\in \calp_{\r}(D)$.

\begin{defi}
\label{p-mes}
A $\calp_{\r}$-measure of a subset $E$ of a $\calp_{\r}$-hyperconvex $\r$-open set $D$ is the upper regularization 
$\omega^*$ of the function 
\begin{equation}
\eqqno(p-mes1)
\omega (z,E,D) = \sup\{u(z):u\in \calp_{\r}^-(D), \text{ } u|_E\le -1\},
\end{equation}
\ie $\omega^* (z,E,D) \deff \limsup_{\xi \to z}\omega (\xi ,E,D)$.
\end{defi}
Remark that $\omega^*$ is $\r$-psh in $D$.  
\begin{defi}
$E$ is called globally $\calp_{\r}$-polar if there exist $u\in \calp_{\r}(D)$, 
which is not $\equiv -\infty$, such that $E\subset \{z\in D:u(z)=-\infty\}$.
\end{defi}

\begin{exmp} 
\label{d-pol-set}
For any $\eta >1$ polydisk $\bar\Delta^{\infty}\left(0, \frac{1}{\n^{1+\eta}}\right)$ is a pluripolar subset of 
$\Delta^{\infty}\left(0, \frac{\adyn}{\n}\right)$. \rm Let's start with proving that function $v$ defined  for 
$z\in \bar\Delta^{\infty}\left(z^0, \r\right)$ as
\begin{equation}
 \eqqno(d-psh-f2)
 v(z) \deff \ln\left(\limsup_{\substack{n\to\infty}}\left\{\frac{|z_n-z^0_n|}{r_n}\right\}\right)
\end{equation}
is $\r$-psh. Indeed, $v$ is a decreasing limit of the $\r$-psh, by the Example \ref{reg-set1}, functions
\[
 v_N(z) \deff \ln\left(\sup_{\substack{n\ge N}}\left\{\frac{|z_n-z^0_n|}{r_n}\right\}\right).
\]
And therefore is $\r$-psh itself. Now remark that for $z\in \bar\Delta^{\infty}\left(0, \frac{1}{\n^{1+\eta}}\right)$
one has that $e^{v_N(z)}\le \frac{r_n^{1+\eta}}{r_n}\to 0$ as $N\to \infty$. Therefore $v_N(z)\searrow -\infty$ for such $z$.
I.e., $v|_{\bar\Delta^{\infty}\left(0, \frac{1}{\n^{1+\eta}}\right)} = -\infty$. At the same time for $z_{\beta}=\{\beta
r_n\}_{n=1}^{\infty}$ one has that $v(z_{\beta}) = \ln \beta\not=-\infty$, \ie $v$ is not $\equiv -\infty$.
\end{exmp}

\begin{defi} 
\label{p-reg-s}
Let $E$ be a subset of a $\calp_{\r}$-hyperconvex open set $D$. $E$ is called locally $\calp_{\r}$-regular if for every 
accumulation  point $z^0$ of $E$ in $D$ and every $\beta>0$ such that  $B^{\infty}_{\r}(z^0,\beta)\subset D$ one has 
\begin{equation}
\omega^*(z^0, E\cap B^{\infty}_{\r}(z^0,\beta), B^{\infty}_{\r}(z^0,\beta)) = -1.
\end{equation}
\end{defi}

\begin{exmp} \rm 
(a) The set (real part) $\el_{r}^{\rr}\cap B^{\infty}_{\r}(0, \beta)$ is locally $\calp_{\r}$-regular in 
$B^{\infty}_{\r}(0, \beta)$. 

\smallskip\noindent (b) For $0<\beta_2<\beta_1<\infty$ the congruent ball $B^{\infty}_{\r}(z^0,\beta_2)$ is a 
$\calp_{\r}$-regular subset of 

\quad the ball $B^{\infty}_{\r}(0,\beta_1)$ provided $\norm{z^0}_{\r} < \beta_1-\beta_2$.
\end{exmp}
The following statement is the so called Two Constants Theorem.

\begin{prop}
 \label{2-c-t}
 Let $D$ be an $\r$-open set in $\el_{\r}$ and let $u\in \calp_{\r}(D)$ be bounded by $M$ on the whole of $D$. 
 Suppose that  for  a subset $E\subset D$ one has that $u|_E\le m$. Then for all $z\in $ the following unequality 
 holds
 \begin{equation}
  \eqqno(2-c-e)
  u(z) \le M\left[1 + \omega^*(z,E,D)\right] - m\omega^*(z,E,D).
\end{equation}
\end{prop}
For the proof remark that function 
\[
v(z) \deff \frac{u(z) -M}{M-m} \quad\text{ belongs to } \calp_{\r}^-(D) \quad\text{ and }\quad v|_E\le -1.
\]
Therefore $v(z) \le \omega^*(z,E,D)$ for all $z\in D$, and this implies \eqqref(2-c-e).

\begin{defi} Call a subset $E\subset D$ of an $\r$-open set $\calp_{\r}$-thick (plurithick) if for any sequence 
$\{u_k\}$ of uniformly bounded from above $\r$-psh functions on $D$ 
\begin{equation}
\eqqno(p-thick-2)
\limsup_{\substack{k\to \infty}}u_k|_E \equiv - \infty \quad\Rightarrow\quad 
\limsup_{\substack{k\to \infty}}u_k \equiv - \infty .
\end{equation}
\end{defi}

\begin{exmp} 
\label{p-thick-3}
(a) For every $z^0\in B^{\infty}_{\r}(0, \beta)$ and every $0<\beta_1<\beta - \norm{z^0}_{\r}$ the set 
$B^{\infty}_{\r}(z^0, \beta_1)$ is $\calp_{\r}$-thick in $B^{\infty}_{\r}(0, \beta)$. \rm By the pathwise
connectivity of the ball it is sufficient to prove this when $z^0=0$, \ie that 
$B^{\infty}_{\r}(0, \beta_1)$ is $\calp_{\r}$-thick in $B^{\infty}_{\r}(0, \beta)$ for $0<\beta_1<\beta$. 
Take any $z\in \bar B^{\infty}_{\r}(0, \beta)$ with $\norm{z}_{\r}=\beta$ and consider the disk 
$\Delta = \{tz\in B^{\infty}_{\r}(0, \beta)\} = i(\{t:|t |< 1\})$, where $i : t\to tz$ is the 
obvious imbedding of $\Delta$ to $B^{\infty}_{\r}(0, \beta)$. Disk $\Delta$ contains a subdisk  
$\Delta_{\frac{\beta_1}{\beta}} \deff i(\{t:|t|<\frac{\beta_1}{\beta}\})=\{tz\in B^{\infty}_{\r}(0,\beta_1)\}$,
which is a thick subset of $\Delta$ in the calssical (\ie one variable) sense. In particular, since 
$\limsup_{\substack{k\to \infty}} u_k|_{\Delta_{\frac{\beta_1}{\beta}}} \equiv - \infty$, it follows that 
$\limsup_{\substack{k\to \infty}}u_k|_{\Delta} \equiv - \infty $, and we conclude that 
$\limsup_{\substack{k\to \infty}}u_k(z) \equiv - \infty $.

\smallskip\noindent (b) The set $\rr^k\cap B^k(0,\beta)$ is plurithick in 
$B^k(0,\beta)$. For $k=1$ this follows from the Poisson formula for harmonic functions. 
The general case canbe easily proved by induction on $k$.
\end{exmp}

\begin{lem}
\label{siciak}
Let $f$ be a holomorphic function on $B^{\infty}_{\r}(z^0,\beta)\times \Delta$. Suppose that for every $z$ in 
some $\calp_{\r}$-thick subset $E$ of $B^{\infty}_{\r}(z^0,\beta)$ the restriction $f_z(s) \deff f(z,s)$ extends
holomorphically to $\cc$. Then $f$ extends to a holomorphic function on $B^{\infty}_{\r}(z^0,\beta)\times\cc$.
\end{lem}
\proof Write the Taylor expansion of $f$ in $B^{\infty}_{\r}(z^0,1)\times \Delta$ in the form
\begin{equation}
\eqqno(taylor4)
f(z^0+z,s) = f(z^0,0) + \sum_{m=1}^{\infty}a_k(z)s^k,
\end{equation}
where $a_k$ are $\r$-holomorphic in $B^{\infty}_{\r}(0,1)$. To get this from \eqqref(taylor1) represent vectors $w$
from $\el_{\r}\times \cc$ as $w = z + se_0$, where $e_0$ spans the subspace $\{0\}\times \cc$ of $\el_{\r}\times \cc$.
With $w^0 = z^0 + 0e_0 = z^0$ write \eqqref(taylor1) in coordinates $w$ as 
\[
 f(w^0 + w) -f(w^0) = \sum_{m=1}^{\infty}P_m(w) = \sum_{m=1}^{\infty}\hat P_m(\underbrace{w,...,w}_{m-\text{times}}) = 
 \sum_{m=1}^{\infty}\hat P_m(\underbrace{z +se_0,...,z + se_0}_{m-\text{times}}) = 
\]
\begin{equation}
 \eqqno(taylor5)
= \sum_{m=1}^{\infty}\sum_{k=0}^{m}C_m^k\hat P_m(\underbrace{se_0,...,se_0}_{k-\text{times}},\underbrace{z,...,z}_{(m-k)-\text{times}})
= \sum_{k=0}^{\infty}\left(\sum_{m=k}^{\infty} C_m^k\hat P_m(\underbrace{e_0,...,e_0}_{k-\text{times}},\underbrace{z,...,z}_{(m-k)-\text{times}})
\right)s^k,
\end{equation}
which gives us \eqqref(taylor4). To estimate the coefficients 
\begin{equation}
 \eqqno(taylor-ak)
 a_k(z) = \sum_{m=k}^{\infty} C_m^k\hat P_m(\underbrace{e_0,...,e_0}_{k-\text{times}},\underbrace{z,...,z}_{(m-k)-\text{times}})
\end{equation}
and justify the change of summution in \eqqref(taylor5)  remark that from \eqqref(a-m-bound) we get that
\[
|\hat P_m(w,...,w)|\le \left(\frac{\beta_1}{\beta_2}\right)^mM_{\beta_2}(f) \quad\text{ on }\quad \Delta^{\infty}(0,\beta_1\r)\times 
\Delta (\beta_1)
\]
for every choice of $0<\beta_1<\beta_2<1$. This implies the uniform on  $\Delta^{\infty}(0,\beta_1\r)$ estimate 
\begin{equation}
 \eqqno(a-est1)
 |a_k(z)| \le \frac{1}{\beta_1^k}\sum_{m=k}^{\infty}\left(\frac{\beta_1}{\beta_2}\right)^{m}M_{\beta_2} =
 \beta_2^k\sum_{m=k}^{\infty}\left(\frac{\beta_1}{\beta_2}\right)^{m-k}M_{\beta_2},
\end{equation}
which implies the normal  convergence of the series \eqqref(taylor5). Therefore
\begin{equation}
\eqqno(a-est2)
\ln \root \k \of{|a_k(z)|} \le 1/k\ln \left(\frac{\beta_2}{(1-\beta_2)(\beta_2-\beta_1)}M_{\beta_2}\right)
\quad\text{ on } \quad \Delta^{\infty}(0,\beta_1\r).
\end{equation}
Functions $\ln \root \k \of{|a_k(z)|}$ are $\r$-psh on $B^{\infty}_{\r}(0,1)$ and tend to $-\infty$ on $E$. 
Due to the assumption of our lemma they tend to $-\infty$ also on $\Delta^{\infty}(0,\beta_1\r)$. We conclude that series 
\eqqref(taylor4) 
converge on $\cc$ for every $z\in B^{\infty}_{\r}(0,1)$. This gives us an extension of $f$ to $z\in B^{\infty}_{\r}(z_0,1)
\times \cc$.

\smallskip To see that this extension (denoted as $\tilde f$) is G\^ateaux differentiable take any two vectors $w^0=z^0+s^0e_0,
w=z+se_0\in B^{\infty}_{\r}(0,1)\times \cc$, take $0<\beta <1$ such that $z^0,z\in B^{\infty}_{\r}(0,\beta )$ and consider 
the two (or one)-dimensional subspace $L$ of $\el_{\r}$ spanned by $z^0$ and $z$. For $t$ in some neighborhood of zero 
vectors $z^0 + tz$ belong to $L\cap B^{\infty}_{\r}(0, 1)$. Applying the calssical Hartogs theorem to the restriction of 
$\tilde f$ to $\left(L\cap B^{\infty}_{\r}(0, 1)\right)\times \cc$ we get the holomorphicity of this restriction. In 
particular $\tilde f(z^0 + tz)$ is holomorphic in $t$. G\^ateaux differentiability is proved.

\smallskip Our extension $\tilde f$ is G\^ateaux differentiable on $B^{\infty}_{\r}(z_0,1)\times \cc$ and holomorphic in 
$B^{\infty}_{\r}(z_0, 1)$ $\times \Delta$, therefore by Theorem 36.5 from \cite{Mu} $\tilde f$ is holomorphic everywhere.
Lemma is proved.

\smallskip\qed
We need also the following form of the previous result.

\begin{corol}
\label{siciak-c}
Let $f$ be a holomorphic function on $B^{\infty}_{\r}(z^0,\beta)\times H_0$. Suppose that for every $z$ in 
some $\calp_{\r}$-thick subset $E$ of $B^{\infty}_{\r}(z^0,\beta)$ the restriction $f_z(s) \deff f(z,s)$ extends
holomorphically to $\cc$. Then $f$ extends to a holomorphic function on $B^{\infty}_{\r}(z^0,\beta)\times\cc$.
\end{corol}
\proof Fractional linear trasformation $\phi (s) = \frac{s-1}{s+1}$ sends $H_0$ to $\Delta$ and $\infty$ to $\adyn$. 
Composing $f$ with $\phi^{-1}$ we get a function $g$ that for every $z\in E$ extends from $\Delta$ to $\cc\setminus 
\{\adyn\}$. Let $\pi :\cc\to \cc\setminus \{\adyn\}$ be the universal covering map. Then $\pi^{-1}\circ g$
satisfies the assumptions of Lemma \ref{siciak}. Therefore it extends to $B^{\infty}_{\r}(z^0,\beta)\times\cc$.
And thus $g$ extends to $B^{\infty}_{\r}(z^0,\beta)\times (\cc\setminus \{\adyn\})$. Which gives the desired
extension of $f$.

\smallskip\qed

\smallskip\noindent{\slsf Proof of Theorem \ref{hart-sic}.} Now let us give the proof of this theorem. 
For every $\z\in B_{\r_2}(0,1)$ we can apply Lemma \ref{siciak} to extend our $f$ to $B_{\r_1}\times \<\cc\cdot \z\>$.
This gives us an extension $\tilde f$ of $f$  to  $B_{\r_1}\times \el_{\r_2}$. In view of already cited Theorem 36.5 
from \cite{Mu} all we need to prove is the G\^ateaux differentiability of this extension. For this take 
any pair of vectors $\z_0,\z_1\in B_{\r_1}(0,1)$ and any $\w_0, \w_1 \in \el_{\r_2}$. Consider the subspace $L$ 
of $\el_{r_1}\times \el_{\r_2}$ generated by $\z_0, \z_1, \w_0,\w_1$. The classical Hartogs theorem gives us 
holomorphicity of $\tilde f|_L$. And therefore we obtain differentiability of $\tilde f(\z_0 +t\z_1, \w_0 + \tau \w_1)$ 
with respect to $t, \tau$. Theorem is proved.

\ifx\undefined\bysame
\newcommand{\bysame}{\leavevmode\hbox to3em{\hrulefill}\,}
\fi

\def\entry#1#2#3#4\par{\bibitem[#1]{#1}
{\textsc{#2 }}{\sl{#3} }#4\par\vskip2pt}
%{ref}{author}{title}ref.

\end{document}